\documentclass[12pt,reqno]{amsart}
\usepackage{amssymb}
\usepackage{amsxtra}
\usepackage{mathrsfs}
\usepackage[all]{xy}
\usepackage[hypertex]{hyperref}
\newtheorem{theorem}{Theorem}[section]
\newtheorem{proposition}[theorem]{Proposition}
\newtheorem{lemma}[theorem]{Lemma}
\newtheorem{corollary}[theorem]{Corollary}
\newtheorem{ex}[theorem]{Example}

\theoremstyle{definition}
\newtheorem{definition}[theorem]{Definition}
\newtheorem{example}[theorem]{Example}
\newtheorem{remark}[theorem]{Remark}

\newcommand*{\ptens}[1]{\mathop{\widehat\otimes}\limits\sb {#1}}
\newcommand*{\Tens}{\mathop{\otimes}}
\newcommand*{\lan}{\mathop{\mathrm{lan}}}
\newcommand*{\ran}{\mathop{\mathrm{ran}}}
\newcommand*{\N}{\mathbb N}
\newcommand*{\B}{\mathcal B}
\newcommand*{\CC}{\mathbb C}
\newcommand*{\Ptens}{\mathop{\widehat\otimes}}
\newcommand*{\id}{\mathbf 1}
\renewcommand*{\Im}{\mathop{\mathrm{Im}}}

\DeclareMathOperator{\Ker}{Ker} 

\begin{document}



\title[Homologically trivial locally $C\sp{\ast}$-algebras]{Structure theory
of homologically trivial \\ and
annihilator locally $C\sp{\ast}$-algebras}

\author{Alexei Yu. Pirkovskii}
\address{Department of Nonlinear Analysis and Optimization\\
Faculty of Science, Peoples' Friendship University of Russia\\
Mikluho-Maklaya 6, Moscow 117198, Russia\\
E-mail: pirkosha@sci.pfu.edu.ru, pirkosha@online.ru}

\author{Yurii V. Selivanov}
\address{Department of Higher Mathematics\\ Russian State
Technological University (MATI)\\ Orshanskaya 3, Moscow 121552, Russia\\
E-mail: selkatya@post.ru}

\dedicatory{Dedicated to the memory of M.~A.~Na\u\i mark on the
centenary of his birth}

\maketitle

\begin{abstract}
We study the structure of certain classes of homologically trivial
locally $C\sp{\ast}$-algebras. These include algebras with
projective irreducible Hermitian $A$-modules, biprojective
algebras, and superbiprojective algebras. We prove that, if $A$ is
a locally $C\sp{\ast}$-algebra, then all irreducible Hermitian
$A$-modules are projective if and only if $A$ is a direct
topological sum of elementary $C\sp{\ast}$-algebras. This is also
equivalent to $A$ being an annihilator (dual, complemented, left
quasi-complemented, or topologically modular annihilator)
topological algebra. We characterize all annihilator
$\sigma$-$C\sp{\ast}$-algebras and describe the structure of
biprojective locally $C\sp{\ast}$-algebras. Also, we present an
example of a biprojective locally $C\sp{\ast}$\mbox{-}algebra that
is not topologically isomorphic to a Cartesian product of
biprojective $C\sp{\ast}$-algebras. Finally, we show that every
superbiprojective locally $C\sp{\ast}$-algebra is topologically
${}\sp{\ast}$\mbox{-}isomorphic to a Cartesian product of full
matrix algebras.
\end{abstract}

\section{Introduction}
This paper is devoted to the study of the structure of locally
$C\sp{\ast}$-algebras satisfying various homological triviality
conditions. The properties of projectivity of topological modules
and of biprojectivity or superbiprojectivity of topological
algebras will play here a central role. We are concerned with the
following questions.

\medskip

\emph{How to characterize, in inner terms of their structure,
locally $C\sp{\ast}$-algebras $A$ such that:}
\begin{itemize}
\item[\rm(i)] \emph{all irreducible Hermitian $A$-modules%
\footnote[1]{\thinspace See Definition \ref{definition3.14}.} are
projective?}

\item[\rm(ii)] \emph{all Hermitian $A$-module are
projective?}

\item[\rm(iii)] \emph{the algebra $A$ is biprojective?}

\item[\rm(iv)] \emph{the algebra $A$ is superbiprojective?}

\end{itemize}

Earlier, similar questions were answered for $C\sp{\ast}$-algebras
\cite{Hel72,Sel76,Sel77,Sel79,Hel97,Hel98,Sel2000,Sel2001}. Apart
from this, the latter two questions were answered for
$\sigma$-$C\sp{\ast}$-algebras \cite{PiS2007}.

The answers to the above questions will be given in Theorems
\ref{theorem3.24}, \ref{theorem4.20}, \ref{theorem5.10} below.

First of all we describe  (in Theorem \ref{theorem3.3}) the
structure of locally $C\sp{\ast}$-algebras with dense socle, and
also the structure of their closed two-sided ideals. Next we show
that property (i) characterizes locally $C\sp{\ast}$-algebras with
dense socle or, equivalently, direct
topological sums of elementary%
\footnote[2]{\thinspace Recall that a $C\sp{\ast}$-algebra $A$ is
called \emph{elementary} if $A$ is isomorphic to the algebra
$\mathcal{K}(H)$ of all compact operators on a Hilbert space $H$.}
$C\sp{\ast}$-algebras. Moreover, the same property is equivalent
to $A$ being an annihilator (dual in Kaplansky's sense,
complemented, left quasi-complemented, or topologically modular
annihilator%
\footnote[3]{\thinspace See Definition \ref{definition3.23}.}%
) topological algebra. Also, we characterize (in Corollary
\ref{corollary3.21}) projective irreducible Hermitian modules over
locally $C\sp{\ast}$-algebras, and describe (in Theorem
\ref{theorem3.29} and Corollary \ref{corollary3.31}) annihilator
unital locally $C\sp{\ast}$-algebras and annihilator
$\sigma$-$C\sp{\ast}$-algebras. In particular, we show that every
annihilator $\sigma$-$C\sp{\ast}$-algebra is topologically
${}\sp{\ast}$-isomorphic to the Cartesian product of a countable
family of annihilator $C\sp{\ast}$-algebras.

Properties (ii) and (iii) turn out to be equivalent. We describe
the structure of biprojective locally
$C\sp{\ast}$\mbox{-}algebras, and present an example of a
biprojective locally $C\sp{\ast}$\mbox{-}algebra that is not
topologically isomorphic to a Cartesian product of biprojective
$C\sp{\ast}$-algebras.

Finally, we show that property (iv) characterizes Cartesian
products of full matrix algebras. In particular, we establish that
all superbiprojective locally $C\sp{\ast}$-algebras are
contractible.

\section{Preliminaries}
By a \emph{topological algebra} we shall mean a complete Hausdorff
locally convex space over $\mathbb{C}$ equipped with a jointly
continuous multiplication. No commutativity or existence of an
identity is assumed.

The symbol $\Ptens$ denotes the complete projective topological
tensor product. The \emph{pro\-duct map} is the continuous linear
operator $\pi\sb A\colon A\Ptens A\to A$ uniquely determined by
$\pi\sb A(a\Tens b)=ab$. By a \emph{Fr\'echet algebra} we mean a
metrizable topological algebra.

We recall that a seminorm $p$ on an algebra $A$ is
\emph{submultiplicative} if $p(ab)\le p(a) p(b)$ for all $a,b\in
A$. If $A$ is a ${}\sp{\ast}$-algebra, then a
\emph{$C\sp{\ast}$-seminorm} on $A$ is a submultiplicative
seminorm $p$ satisfying ${p(a\sp{\ast})=p(a)}$ and ${p(a\sp{\ast}
a)=p(a)\sp 2}$ for all $a\in A$. A topological algebra $A$ is an
\emph{Arens-Michael algebra} if the topology on $A$ can be defined
by a family of submultiplicative seminorms.
A \emph{locally $C\sp{\ast}$-algebra}%
\footnote[4]{\thinspace These objects are called
$b\sp{\ast}$-algebras in \cite{Allan,Dixon,Apos},
$LMC\sp{\ast}$-algebras in \cite{Lassner,Schmdgn},
pro\mbox{-}$C\sp{\ast}$\mbox{-}algebras in
\cite{Voic,Phillips1,Phillips2} and multi-$C\sp{\ast}$-algebras in
\cite{Hel93}.} is an Arens-Michael ${}\sp{\ast}$-algebra whose
topology is defined by a family of $C\sp{\ast}$\mbox{-}seminorms.
Recall that the term ``locally $C\sp{\ast}$-algebra'' is due to
A.~Inoue \cite{Inoue}. Metrizable locally $C\sp{\ast}$-algebras
are also called \emph{$\sigma$-$C\sp{\ast}$-algebras}%
\footnote[5]{\thinspace Also, metrizable locally
$C\sp{\ast}$-algebras are called $F\sp{\ast}$-algebras in
\cite{Brooks}.} \cite{Arves,Phillips1,Phillips2}.

As is known, every locally $C\sp{\ast}$\mbox{-}algebra is
topologically ${}\sp{\ast}$-isomorphic to an inverse limit of
$C\sp{\ast}$\mbox{-}algebras. $\sigma$-$C\sp{\ast}$-algebras are
countable inverse limits of $C\sp{\ast}$\mbox{-}algebras.

S. J. Bhatt and D. J. Karia obtained an interesting
characterization of locally $C\sp{\ast}$-algebras among
Arens-Michael ${}\sp{\ast}$-algebras:

\begin{theorem}[{see \cite{Bhatt-Karia}}]
An Arens-Michael ${}\sp{\ast}$-algebra $A$ is a locally
$C\sp{\ast}$-algebra if and only if there exists a dense
${}\sp{\ast}$-subalgebra $B$ of $A$ which is a
$C\sp{\ast}$-algebra under some norm and which is continuously
embedded in $A$.
\end{theorem}

We recall (cf.~\cite[Lemma 8.14]{Frag1}) that every locally
$C\sp{\ast}$-algebra $A$ is semisimple, i.e., the Jacobson radical
of $A$ is~$0$. For details and references on locally
$C\sp{\ast}$-algebras, see \cite{Frag3,Phillips1,Phillips2}.

\medskip

Now for any subset $S$ of an algebra $A$, let $\lan(S)$ and
$\ran(S)$ denote the \emph{left} and \emph{right annihilators} of
$S$ in $A$, respectively (see, e.g., \cite[\S30]{Bon}). So we have
\[
\lan(S)=\{a\in A\mid ab=0\mbox{ for all }b\in
S\}\quad\mbox{and}\quad\ran(S)=\{a\in A\mid ba=0\hbox{ for all
}b\in S\}.
\]
As is known, $\lan(S)$ is always a left ideal of $A$, and
$\ran(S)$ is always a right ideal of~$A$. Furthermore, the left
(respectively, right)
annihilator of a left (respectively, right) ideal
is a two-sided ideal. Also, if $A$ is a topological algebra,
then both $\lan(S)$ and $\ran(S)$ are closed.

Recall that a net $\{e\sb {\nu},\,\nu\in \Lambda\}$ in a
topological algebra $A$ is called a \emph{left} (respectively,
\emph{right}) \emph{bounded approximate identity} if $e\sb {\nu}a
\to a$ (respectively, $ae\sb {\nu} \to a$) for each $a\in A$ and
if the elements $e\sb {\nu}$ form a bounded set in $A$.
A \emph{bounded approximate identity} is a net which is both
a left and a right bounded approximate identity.

We recall that a locally $C\sp{\ast}$-algebra always has a bounded
approximate identity (see \cite[Theorem 11.5]{Frag3}). This
implies, in particular, that $\lan(A)=\ran(A)=0$ for every such
algebra $A$.

A topological algebra $A$ is an \emph{annihilator algebra} if,
for every closed left ideal $J$ and
for every closed right ideal $K$, we have $\ran(J)=0$ if and only
if $J=A$ and $\lan(K)=0$ if and only if $K=A$. If
$\lan(\ran(J))=J$ and $\ran(\lan(K))=K$, then $A$ is called a
\emph{dual algebra}. It is obvious that a dual algebra is
automatically an annihilator algebra.

Dual algebras were introduced by I.~Kaplansky \cite{Kapl48}, and
annihilator algebras were introduced by F.~F.~Bonsall and
A.~W.~Goldie \cite{BoGo}. The existence of an annihilator Banach
algebra which is not dual was first established by B.~E.~Johnson
\cite{Joh67}. His example was commutative and semisimple. Later A.
M. Davie \cite{Davie} gave an example of a topologically simple
annihilator Banach algebra which is not dual. However every
annihilator $C\sp{\ast}$-algebra is dual (see \cite[Corollary
4.10.26]{Rick}).

Now we recall (see \cite{Sel2000} or \cite[p.~63]{Hel2006}) that
the $c\sb 0$-{\it sum}, $\bigoplus\nolimits\sb 0\{A\sb \nu\mid
\nu\in \Lambda\}$, of a family of $C\sp{\ast}$-algebras $A\sb
\nu$, $\nu\in \Lambda$, is defined to be the set of all functions
$f$ defined on $\Lambda$ such that
\begin{itemize}
\item[\rm(i)] $f(\nu)\in A\sb \nu$ for each $\nu\in \Lambda$, and

\item[\rm(ii)] for each
$\varepsilon>0$ the set $\{\nu\,:\,\|f(\nu)\|\ge\varepsilon\}$ is
finite.
\end{itemize}
As is known, this set is a $C\sp{\ast}$-algebra with respect to
pointwise operations and the norm $\|f\|=\sup\limits\sb
{\nu\in\Lambda}\|f(\nu)\|$.

Following \cite{Dixm}, we say that a $C\sp{\ast}$-algebra $A$ is
\emph{elementary} if $A$ is isometrically ${}\sp{\ast}$-isomorphic
to the algebra $\mathcal{K}(H)$ of all compact operators on a
Hilbert space $H$. We recall (see \cite[4.7.20]{Dixm}) that a
$C\sp{\ast}$-algebra $A$ is annihilator (or, equivalently, dual)
if and only if it is isometrically ${}\sp{\ast}$-isomorphic to the
$c\sb 0$-sum of a family of elementary $C\sp{\ast}$-algebras.

For details and references on annihilator and dual algebras, see
\cite{Rick,Naim,Bon,Pal1}.

\medskip

Let $A$ be a topological algebra and let $\mathcal{M}\sb\ell$ be
the set of all closed left ideals of~$A$. Then $A$ is called a
\emph{left quasi-complemented algebra} if there exists a mapping
$q\colon J\mapsto J\sp q$ of $\mathcal{M}\sb\ell$ into itself
having the following properties:
\begin{itemize}
\item[\rm(i)] $J\cap J\sp q=0\ (J\in \mathcal{M}\sb\ell)$;

\item[\rm(ii)] $(J\sp q)\sp q=J\ (J\in \mathcal{M}\sb\ell)$;

\item[\rm(iii)] if $J\sb 1\subset J\sb 2$, then
$J\sb 2\sp q\subset J\sb 1\sp q\ (J\sb 1,J\sb 2\in
\mathcal{M}\sb\ell)$.
\end{itemize}
The mapping $q$ is called a \emph{left quasi-complementor} on $A$.

A left quasi-complemented algebra is called a \emph{left
complemented algebra} if it satisfies:
\begin{itemize}
\item[\rm(iv)] $J+J\sp q=A\ (J\in \mathcal{M}\sb\ell)$.
\end{itemize}
In this case, the mapping $q$ is called a \emph{left complementor}
on $A$.

\emph{Right quasi-complemented algebras} and \emph{right
complemented algebras} are defined analogously. A left and right
complemented algebra is called a \emph{complemented algebra}.

Complemented Banach algebras were introduced by B.~J.~Tomiuk
\cite{Tom62} and have been studied by many authors (see, e.g.,
\cite{Olub,F_Alex1,Al-Tom69,F_Alex2,Wong70,Tom_Wo,F_Alex3,Wong74,Malviya,Tom88,Tom92}).
Left (right) quasi-complemented topological algebras were defined
by T.~Husain and P.~K.~Wong \cite{HuWo72}. They showed that there
exist right quasi-complemented algebras which are not right
complemented.

Structural properties of left complemented semisimple topological
algebras in which every modular maximal left ideal is closed were
studied by M. Haralampidou \cite{Har93}.

It is known that a $C\sp{\ast}$-algebra is complemented if and
only if it is dual (see \cite{Al-Tom69}). A similar result is true
for left (right) quasi-complemented algebras (see \cite{HuWo72}).

Let $\{I\sb {\nu}\mid \nu\in \Lambda\}$ be a family of (left,
two-sided) ideals in an algebra $A$. Recall that the smallest
(left, two-sided) ideal in $A$ which contains every $I\sb {\nu}$
is called the \emph{sum} of the ideals $I\sb {\nu}$. The sum of
the ideals $I\sb {\nu}$ evidently consists of all finite sums of
elements from the ideals $I\sb {\nu}$.

If $A$ is a topological algebra, then the closure of the sum of
the ideals $I\sb {\nu}$ is called their \emph{topological sum}. If
each $I\sb {\nu}$ is closed and intersects the topological sum of
the remaining ideals in the zero element, then the topological sum
is called a \emph{direct topological sum} (C.~E.~Rickart's terminology,
see \cite[p.~46]{Rick}).

We recall that a \emph{minimal left ideal} of an algebra $A$ is a
left ideal $J\neq 0$ such that $0$~and $J$ are the only left
ideals contained in $J$. Recall also (see, e.g., \cite[\S30]{Bon})
that the \emph{left socle} of an algebra $A$ is the sum of all
minimal left ideals of $A$. The \emph{right socle} is similarly
defined in terms of right ideals. If the left socle coincides with
the right socle, then it is called the \emph{socle} of $A$ and is
denoted by $\mathrm{Soc}(A)$.

Recall that an algebra $A$ is \emph{semiprime} if it has no
two-sided ideals, $I\neq 0$, with $I\sp 2=0$. Since semisimple
algebras are semiprime (see, e.g., \cite[Proposition~30.5]{Bon}),
we conclude that every locally $C\sp{\ast}$-algebra is semiprime.

Recall that, if $A$ is semiprime, then the left and right socles
of $A$ coincide, and so $A$ has a~socle. $\mathrm{Soc}(A)$ is
known to be a two-sided ideal of $A$. Recall that a
$C\sp{\ast}$-algebra $A$ is dual if and only if $\mathrm{Soc}(A)$
is dense in $A$ (see \cite[4.7.20]{Dixm}).

Also, we recall that a \emph{minimal closed two-sided ideal} of a
topological algebra is a closed two-sided ideal $I\neq 0$ that
contains no closed two-sided ideals other than $0$ and~$I$.
Finally, we recall that a topological algebra $A$ is
\emph{topologically simple} if its only closed two-sided ideals
are $0$ and $A$.

The proof of the following lemma repeats the proof of
\cite[Lemma~32.4]{Bon}.

\begin{lemma}
Let $A$ be a semiprime topological algebra, and let $I$ be a
two-sided ideal of~$A$. Then $\lan (I)=\ran (I)$, $I\cap\lan
(I)=0$. If, in addition, $A$ is an annihilator algebra, then
$I\oplus\lan (I)$ is dense in $A$.
\end{lemma}

If $A$ is a locally $C\sp{\ast}$-algebra, then even more is true.

\begin{lemma}[{cf.~\cite[\S25.2]{Naim} and
\cite{Har2005}}]\label{lemma2.3} Let $A$ be a locally
$C\sp{\ast}$-algebra.
\begin{itemize}
\item[\rm(i)] If $J$ is a closed left ideal of $A$, then
$J\cap{}(\ran (J))\sp{\ast}=0$ and $J\oplus {}(\ran (J))\sp{\ast}$
is a closed left ideal of $A$. If, in addition, $A$ is an
annihilator algebra, then $A=J\oplus{}(\ran (J))\sp{\ast}$.
\item[\rm(ii)] If $K$ is a closed right ideal of $A$, then
$K\cap {}(\lan (K))\sp{\ast}=0$ and $K\oplus {}(\lan
(K))\sp{\ast}$ is a closed right ideal of $A$. If, in addition,
$A$ is an annihilator algebra, then $A=K\oplus {}(\lan
(K))\sp{\ast}$.
\item[\rm(iii)] If $I$ is a closed two-sided ideal of $A$, then
$I\cap\lan (I)=0$ and $I\oplus \lan (I)$ is a closed two-sided
ideal of $A$. If, in addition, $A$ is an annihilator algebra, then
$A=I\oplus\lan (I)$.
\end{itemize}
\end{lemma}

\begin{proof} (i) Clearly, $(\ran (J))\sp{\ast}$ is a closed
left ideal of $A$. Suppose that $a\in J\cap{}(\ran (J))\sp{\ast}$.
Then $aa\sp{\ast}=0$, and consequently, for each continuous
$C\sp{\ast}$-seminorm $p$ on $A$,
\[
p(a)=p(aa\sp{\ast})\sp {1/2}=0.
\]
Hence $a=0$, and so $J\cap{}(\ran (J))\sp{\ast}=0$.

Now let $L=J\oplus {}(\ran (J))\sp{\ast}$, and let $a=b+c\in L$,
where $b\in J$, $c\in (\ran (J))\sp{\ast}$. Then
$ac\sp{\ast}=cc\sp{\ast}$, and therefore
\[p(c)\sp 2=p(cc\sp{\ast})\le p(a)p(c\sp{\ast})\]
and
\begin{equation}
\label{pcltpa} p(c)\le p(a),
\end{equation}
for each continuous $C\sp{\ast}$-seminorm $p$ on $A$.

Similarly, we have
\begin{equation}
\label{pbltpa} p(b)\le p(a),
\end{equation}
for each continuous $C\sp{\ast}$-seminorm $p$ on $A$. It follows
easily from these inequalities that the left ideal $L$ is closed.

We prove now that $\ran(L)=0$. Indeed, suppose that $a\in
\ran(L)$. Then $Ja=0$ and $(\ran (J))\sp{\ast} a=0$. It follows
from the first equality that $a\in \ran (J)$ and $a\sp{\ast}\in
(\ran (J))\sp{\ast}$. Hence $a\sp{\ast} a=0$, and therefore, for
each continuous $C\sp{\ast}$-seminorm $p$ on $A$,
\[p(a)=p(a\sp{\ast} a)\sp {1/2}=0,\] and so $a=0$. Thus
$\ran(L)=0$, and we see that, if $A$ is an annihilator algebra,
then $L=A$.

(ii) This is similar.

(iii) This follows from (ii) and from the fact that every closed
two-sided ideal of a locally $C\sp{\ast}$-algebra is a
${}\sp{\ast}$-ideal (see \cite[Theorem 11.7]{Frag3}).
\end{proof}

The proof of the next lemma repeats the proof of \cite[Theorem
2.8.29]{Rick}; see also \cite[Theorem 2.4]{Har2005}.

\begin{lemma}\label{lemma2.4}
Let $A$ be a topological algebra such that $\lan(A)=\ran(A)=0$.
Suppose that $A$ is equal to the topological sum of a given family
$\{I\sb {\nu}\mid \nu\in \Lambda\}$ of its closed two-sided
ideals. If each $I\sb {\nu}$ is an annihilator algebra, then $A$
is an annihilator algebra.
\end{lemma}

\section{Locally $C\sp{\ast}$-algebras with projective irreducible
modules}

\subsection{Locally $C\sp{\ast}$-algebras with dense socle}
Let $A$ be a locally $C\sp{\ast}$-algebra, and let $P$ be the
family of all continuous $C\sp{\ast}$-seminorms on $A$. For each
$p\in P$ we set
\[
N\sb p=\Ker p,
\]
where $\Ker p=\{a\in A \mid p(a)=0\}$.  Then each $N\sb p$ is a
two-sided ${}\sp{\ast}$-ideal of $A$. The quotient seminorm of $p$
on $A/N\sb p$ is a $C\sp{\ast}$-norm. We denote this norm by
$\|\cdot\|\sb p$. The completion of $A/N\sb p$ with respect to
this norm is a $C\sp{\ast}$-algebra which is denoted by $A\sb p$.
Following \cite{Hel93}, we call the $C\sp{\ast}$-algebra $A\sb p$
\emph{concomitant} with $A$.

Note that every quotient algebra $A/N\sb p$ is in fact already
complete with respect to the norm $\|\cdot\|\sb p$ (see
\cite{Apos}, \cite{Schmdgn}, \cite[Corollary 1.2.8]{Phillips2}),
and so $A\sb p=(A/N\sb p,\,\|\cdot\|\sb p)$.

\begin{proposition}
\label{proposition3.1} Let $A$ be a locally $C\sp{\ast}$-algebra,
and let $P$ be the family of all continuous $C\sp{\ast}$-seminorms
on $A$. Then, for every $p\in P$, the quotient topology on $A/N\sb
p$ coincides with the topology determined by the norm
$\|\cdot\|\sb p$.
\end{proposition}

\begin{proof}
Evidently, the first topology is stronger than the second one. On
the other hand, let $q\in P$, and suppose that $q\ge p$ (i.e.,
$q(a)\ge p(a)$ for all $a\in A$). Then the quotient seminorm of
$q$ on $A/N\sb p$ is equal to the quotient norm of the
$C\sp{\ast}$-norm $\|\cdot\|\sb q$ on $A\sb q=A/N\sb q$, and hence
is itself a (complete) $C\sp{\ast}$-norm. Therefore the latter
norm is equal to $\|\cdot\|\sb p$.

Note that $P$ is directed with the order $p\le q$ (see
\cite{Phillips2}). Hence we conclude that, if $q\in P$ is
arbitrary, and $\widehat{q}\,$ is the quotient seminorm of $q$ on
$A/N\sb p$, then $\|\cdot\|\sb p\ge \widehat{q}$. The rest is
clear.
\end{proof}

Let $A$ be an algebra. We recall \cite{Bon} that a non-zero
idempotent $e\in A$ is called \emph{minimal} if $eAe$ is a
division algebra. Since an Arens-Michael division algebra
is isomorphic to $\CC$ (see, e.g., \cite[Theorem V.1.7]{Hel93}),
we conclude that, for a minimal idempotent $e$ in an Arens-Michael algebra $A$,
we have $eAe=\CC e$.

\begin{corollary}\label{corollary3.2}
Let $A$ be a topologically simple locally $C\sp{\ast}$-algebra
that contains a mi\-ni\-mal idempotent. Then $A$ is isomorphic, as
a topological ${}\sp{\ast}$-algebra, to the $C\sp{\ast}$-algebra
of all compact operators on a Hilbert space.
\end{corollary}

\begin{proof} Let $p\in P$, and suppose that $p\neq 0$. Then
$\Ker p$ is a closed two-sided ideal of $A$, and hence $\Ker p=0$
because $A$ is topologically simple. Thus $p$ is a norm. By
Proposition~\ref{proposition3.1}, the topology on $A=A/N\sb p$
coincides with the topology determined by the norm $p=\|\cdot\|\sb
p$. So $A$ is isomorphic, as a topological ${}\sp{\ast}$-algebra,
to the $C\sp{\ast}$-algebra $(A,\,\|\cdot\|\sb p)$. Since $A$
contains a minimal idempotent, it follows from
\cite[Proposition~30.6]{Bon} that $\mathrm{Soc}(A)\neq 0$. Since
$A$ is topologically simple, and $\mathrm{Soc}(A)$ is a two-sided
ideal, it follows that $\mathrm{Soc}(A)$ is dense in $A$. It
remains to apply \cite[4.7.20]{Dixm} and \cite[Corollary
4.10.20]{Rick}.
\end{proof}

\begin{theorem}\label{theorem3.3}
Let $A$ be a locally $C\sp{\ast}$-algebra with dense socle, and
let $\{I\sb {\nu}\mid \nu\in \Lambda\}$ be the collection of
minimal closed two-sided ideals of this algebra; we put $P\sb
{\nu}=\lan(I\sb {\nu})$ $(=\ran(I\sb {\nu}))$ $(\nu\in \Lambda)$.
Then the following statements hold.
\begin{itemize}
\item[\rm(i)] $A$ is the direct topological sum of all the $I\sb {\nu}$, and
moreover, for each ${\nu\sb {1},\nu\sb {2}\in \Lambda}$ with
${\nu\sb {1}\neq \nu\sb {2}}$, we have $I\sb {\nu\sb {1}} I\sb {\nu\sb
{2}}=0$.
\item[\rm(ii)] For each $\nu\in \Lambda$ the algebra $I\sb {\nu}$ is
isomorphic, as a topological ${}\sp{\ast}$-algebra, to the
$C\sp{\ast}$-algebra $\mathcal{K}(H\sb {\nu})$, where $H\sb {\nu}$
is a Hilbert space.
\item[\rm(iii)] Every closed two-sided ideal $I$ of $A$ is the intersection of the
ideals $P\sb {\nu}$ that contain~$I$. In particular, $\bigcap
\limits\sb {\nu\in \Lambda}P\sb {\nu}=0$.
\item[\rm(iv)] Every closed two-sided ideal $I$ of $A$ is the direct topological sum
of the ideals $I\sb {\nu}$ contained in $I$. As a consequence, for
each $\nu$, $P\sb {\nu}=\overline{\bigoplus\limits\sb {\mu\in
\Lambda\setminus \{\nu\}} I\sb {\mu}}$.
\item[\rm(v)] $A$ is an annihilator algebra.
\item[\rm(vi)] For each $\nu$, $A=I\sb {\nu}\oplus P\sb {\nu}$.
Moreover, let $\psi\sb {\nu}\colon A\to A/P\sb {\nu}$, $a \mapsto
a\sb {\nu}$, denote the quotient map. Then $\varphi\sb {\nu}\colon
I\sb {\nu}\hookrightarrow A\xrightarrow{\psi\sb {\nu}} A/P\sb
{\nu}$ is an isometric ${}\sp{\ast}$-isomorphism of
$C\sp{\ast}$-algebras.
\item[\rm(vii)] The homomorphism of locally $C\sp{\ast}$-algebras \[
\tau\colon A\longrightarrow \prod\sb {\nu\in \Lambda} A/P\sb
{\nu}=\prod\sb {\nu\in \Lambda} \mathcal{K}(H\sb {\nu}),\quad a
\mapsto \{a\sb {\nu}\},
\]
is a continuous embedding with dense range.
\item[\rm(viii)] If $I$, $I=\overline{\bigoplus\limits\sb {\nu\in \Lambda\sb 1} I\sb
{\nu}}$, is an arbitrary closed two-sided ideal of $A$, then
$A=I\oplus J$, where
\[
J=\lan(I)=\overline{\bigoplus\sb {\nu\in \Lambda\setminus
\Lambda\sb 1} I\sb {\nu}}.
\]
Moreover, the natural homomorphism ${\varphi\colon
J\hookrightarrow A\rightarrow A/I}$ is a topological
${}\sp{\ast}$-iso\-mor\-phism. In particular, $A/I$ is complete.
\item[\rm(ix)] If $I$ is a closed two-sided ideal of $A$, then
$\mathrm{Soc}(I)$ is dense in $I$, $\mathrm{Soc}(A/I)$ is dense
in~$A/I$, and $A$ is topologically ${}\sp{\ast}$-isomorphic to the
Cartesian product of the algebras $I$ and $A/I$.
\item[\rm(x)] If $p$ is a continuous
$C\sp{\ast}$-seminorm on $A$, then there exists a subset
${\Lambda\sb 1\subset\Lambda}$ such that $\Ker p=\lan(J)$, where
$J=\overline{\bigoplus\limits\sb {\nu\in \Lambda\sb 1} I\sb
{\nu}}$ and $J$ is isomorphic, as a topological
${}\sp{\ast}$-algebra, to a $C\sp{\ast}$-algebra (and hence $J$ is
the $c\sb 0$-{\it sum} of the family of $C\sp{\ast}$-algebras
$I\sb {\nu}=\mathcal{K}(H\sb {\nu})$, $\nu\in \Lambda\sb 1$).
Moreover, for each ${a\in A}$, ${p(a)=\sup\limits\sb {\nu\in
\Lambda\sb 1}\left\{\|a\sb {\nu}\|\sb {\mathcal{K}(H\sb
{\nu})}\right\}}$.
\item[\rm(xi)] There exists a ${}\sp{\ast}$-homomorphism of locally $C\sp{\ast}$-algebras
\[\omega\colon \bigoplus\nolimits\sb 0\{\mathcal{K}(H\sb {\nu})\mid
\nu\in \Lambda\} \longrightarrow A \] such that $\omega$ is a
continuous embedding with dense range, and $\tau\circ\omega$ is
the natural embedding of the $C\sp{\ast}$-algebra
$\bigoplus\nolimits\sb 0\{\mathcal{K}(H\sb {\nu})\mid \nu\in
\Lambda\}$ into $\prod\limits\sb {\nu\in \Lambda} \mathcal{K}(H\sb
{\nu})$.

\end{itemize}
\end{theorem}

\begin{proof} (i) Note first that, if $J$ is a minimal left ideal of $A$,
then the closed two-sided ideal $I$ generated by $J$ is a minimal
closed two-sided ideal of $A$.

Indeed (cf.~\cite[\S25.3, II]{Naim}), let $K$ be a closed
two-sided ideal of $A$ contained in $I$. Then $J\cap K$ is a left
ideal of $A$ contained in the minimal left ideal $J$. Therefore
either $J\subset K$ or $J\cap K=0$. In the first case $I\subset
K$, and so $K=I$. In the second case \[K\cdot J\subset J\cap
K=0,\] and consequently $J\subset \ran(K)$. Since $\ran(K)$ is a
closed two-sided ideal of $A$, we have $I\subset \ran(K)$. Hence
$K\subset \ran(K)$, i.e., $K\sp 2=0$. Since $A$ is
semiprime, we get $K=0$. Therefore $I$ is a minimal closed
two-sided ideal of $A$.

Thus $A$ coincides not only with the topological sum
of all its minimal left ideals but also with the topological sum
of all the minimal closed two-sided ideals $I\sb {\nu}$, $\nu\in
\Lambda$. The latter topological sum is direct, and moreover, for
each ${\nu\sb {1},\nu\sb {2}\in \Lambda}$ with ${\nu\sb {1}\neq
\nu\sb {2}}$, we have $I\sb {\nu\sb {1}} I\sb {\nu\sb {2}}=0$
(cf.~\cite[the proof of~Theorem~5 in \S25.3]{Naim}).

(ii) Let $\nu\in \Lambda$. Note first that the minimal closed
two-sided ideal $I\sb {\nu}$ is a topologically simple locally
$C\sp{\ast}$-algebra.

Indeed, $I\sb {\nu}$ is a ${}\sp{\ast}$-ideal by \cite[Theorem
11.7]{Frag3}. Thus $I\sb {\nu}$ is a closed
${}\sp{\ast}$-subalgebra of $A$, i.e., a locally
$C\sp{\ast}$-algebra. Since $A$ is the direct topological sum of
all the $I\sb {\nu}$, it follows easily that $I\sb {\nu}$ is a
topologically simple algebra.

We prove now that $I\sb {\nu}$ contains a minimal left ideal $L$
of $A$. It is then obvious that $I\sb {\nu}$ is generated by $L$.

Indeed, suppose $I\sb {\nu}$ contains no minimal left ideals of
$A$, and let $J$ be a minimal left ideal. Then $J\cap I\sb
{\nu}\subset J$ and $J\cap I\sb {\nu}\neq J$, since $J\not\subset
I\sb {\nu}$. Consequently $J\cap I\sb {\nu}=0$, and therefore
$I\sb {\nu}\cdot J\subset J\cap I\sb {\nu}=0$. Since $J$ is
arbitrary, $I\sb {\nu}\cdot\mathrm{Soc}(A)=0$. Hence
\[
(I\sb {\nu})\sp 2\subset I\sb {\nu}\cdot A= I\sb
{\nu}\cdot\overline{\mathrm{Soc}(A)}=0.
\]
Since $A$ is semiprime, $I\sb {\nu}=0$. So we have a
contradiction.

Finally we note that, by \cite[Lemma 30.2]{Bon}, $L$ has the form
$L=Ae$, where $e\in I\sb {\nu}\subset A$ is a minimal idempotent.
Since $eAe$ is an Arens-Michael division algebra, we have
$eAe=\mathbb Ce$, and hence $eI\sb {\nu}e=\CC e$. It remains to
apply Corollary \ref{corollary3.2}.

(iii) Let $I$ be a closed two-sided ideal of $A$, and let $a\in
\bigcap \{P\sb {\nu}:P\sb {\nu}\supset I\}$. Let $I\sb {\mu}$,
$\mu\in \Lambda$, be a minimal closed two-sided ideal of $A$. Then
either $I\sb {\mu}\subset I$ or $I\sb {\mu}\cap I=0$. In the first
case, $I\sb {\mu}a\subset I\sb {\mu}\subset I$. In the second case,
$I\cdot I\sb {\mu}\subset I\cap I\sb {\mu}=0$, and consequently
$I\subset \lan(I\sb {\mu})=P\sb {\mu}$. Hence $a\in P\sb {\mu}$
and $I\sb {\mu}a\subset I\sb {\mu}\cdot P\sb {\mu}=0\subset I$.
Thus $I\sb {\mu}a\subset I$ for each  $\mu\in \Lambda$. Since the
sum of the ideals $I\sb {\mu}$ is dense in $A$,
$\overline{Aa}\subset I$. Since $A$
has an approximate identity, we see that $a\in I$.

(iv) Let $I$ be a closed two-sided ideal of $A$, and let
$\Lambda\sb 1=\{\nu\in \Lambda:I\sb {\nu}\subset I\}$ and
$J=\overline{\bigoplus\limits\sb {\nu\in \Lambda\sb 1} I\sb
{\nu}}$. As we noted in the proof of (iii), $I\subset P\sb {\nu}$
for each $\nu\in \Lambda\setminus \Lambda\sb 1$, and hence
$J\subset P\sb {\nu}$. For each $\nu\in \Lambda\sb 1$ we have
$J\not\subset P\sb {\nu}$, since otherwise $I\sb {\nu}\subset P\sb
{\nu}$ and $I\sb {\nu}=I\sb {\nu}\cap P\sb {\nu}=0$. By (iii),
$J=\bigcap \limits\sb {\nu\in \Lambda\setminus \Lambda\sb 1} P\sb
{\nu}$. Since $I$ is contained in this intersection, we conclude
that $I=J$.

(v) This follows from Lemma \ref{lemma2.4} and from the fact that
the Banach algebra of all compact operators on a Hilbert space is
an annihilator algebra (see, e.g., \cite[\S25.6, Theorem
13]{Naim}).

(vi) This follows from (v), from Lemma \ref{lemma2.3} and from
inequalities (\ref{pcltpa}) and (\ref{pbltpa}) (see the proof
of Lemma \ref{lemma2.3}).

(vii) This follows easily by using (i)--(iv) and (vi).

(viii) This follows by the same argument as in (vi).

(ix) This follows at once from (viii) and (iv).

(x) Note that $\Ker p$ is a closed two-sided ideal of $A$. By
(iv), there exists a subset $\Lambda\sb 0\subset\Lambda$ such that
$\Ker p=\overline{\bigoplus\limits\sb {\nu\in \Lambda\sb 0}I\sb
{\nu}}$. We set \[\Lambda\sb 1=\Lambda\setminus \Lambda\sb 0.\] By
(viii), we have $A=\Ker p\oplus J$, where \[J=\lan(\Ker
p)=\overline{\bigoplus\limits\sb {\nu\in \Lambda\sb 1} I\sb
{\nu}},\] and moreover the natural homomorphism
\[{\varphi\colon J\hookrightarrow A\rightarrow A/\Ker p}\] is a topological
${}\sp{\ast}$-isomorphism of locally $C\sp{\ast}$-algebras. The
rest follows from the fact that, by
Proposition~\ref{proposition3.1}, $A/\Ker p$ is topologically
${}\sp{\ast}$-isomorphic to the concomitant $C\sp{\ast}$-algebra
$(A\sb p,\,\|\cdot\|\sb p)$.

(xi) This follows easily by using (i), (ii), (vi), (vii) and (x).
\end{proof}

The next corollary follows easily from Theorem \ref{theorem3.3}.

\begin{corollary}\label{corollary3.4}
Let $A$ be a locally $C\sp{\ast}$-algebra with dense socle.
\begin{itemize}
\item[\rm(i)] If, for each continuous $C\sp{\ast}$-seminorm $p$ on $A$, the
concomitant $C\sp{\ast}$-algebra $A\sb p$ has finite spectrum
$\widehat{A}\sb p$, then $A$ is topologically
${}\sp{\ast}$-isomorphic to the Cartesian product of a family of
elementary $C\sp{\ast}$-algebras.
\item[\rm(ii)] If, for each continuous $C\sp{\ast}$-seminorm $p$ on $A$, the
concomitant $C\sp{\ast}$-algebra $A\sb p$ is finite-dimensional,
then $A$ is topologically ${}\sp{\ast}$-isomorphic to the
Cartesian product of a family of full matrix algebras.
\item[\rm(iii)] If $A$ has an identity, then $A$ is topologically
${}\sp{\ast}$-isomorphic to the Cartesian product of a family of
full matrix algebras.
\end{itemize}
\end{corollary}

\subsection{Some definitions and results concerning topological algebras and modules}
Suppose that $A$ is a topological algebra%
\footnote[6]{\thinspace Recall that by a topological algebra we
mean a complete Hausdorff locally convex space equipped with a
jointly continuous multiplication.} and $X$ is a complete
Hausdorff locally convex space which is also a left $A$-module.
Thus there is a bilinear map $(a,x)\mapsto a\cdot x$ from $A\times
X$ to $X$ such that $(ab)\cdot x=a\cdot (b\cdot x)$ for $a,b\in
A$, $x\in X$. Then $X$ is called a \emph{left topological
$A$-module} if the above module map is jointly continuous. For two
such modules $X$ and $Y$, an \emph{$A$-module morphism} from $X$
to $Y$ is a continuous linear operator $\varphi\colon X\to Y$
which is a module homomorphism.

Similar definitions apply to \emph{right topological $A$-modules}
and \emph{topological $A$-bi\-mo\-du\-les}. For example, the
algebra $A$ is itself a topological $A$-bimodule with respect to
the maps given by the product in $A$. If $X$ is a left topological
$A$-module and $Y$ is a right topological $A$-module, then
$X\Ptens Y$ is a topological $A$-bimodule for the products defined
by
\[
a\cdot (x\Tens y)=a\cdot x\Tens y,\quad (x \Tens y)\cdot a=x\Tens
y\cdot a\qquad (a\in A,\: x\in X,\: y\in Y).
\]
In particular, $A\Ptens A$ is a topological $A$-bimodule
in this way.

For a left $A$-module $X$, we denote by $A\cdot X$ the linear span
of the elements of the form $a\cdot x$, where $a\in A$ and $x\in
X$; expressions of the type $Y\cdot A$ have a similar meaning for
right $A$-modules $Y$. We write $A\sp 2$ for $A\cdot A$. A
topological algebra $A$ is \emph{idempotent} if $A\sp
2$ is dense in $A$.

A left topological $A$-module $X$ is
\emph{essential} if $A\cdot X$ is dense in $X$. A left
module $X$ over an algebra $A$ is (algebraically) \emph{irreducible}
if $A\cdot X\neq 0$ and $X$~contains no non-zero proper
submodules. It is obvious that, if $X$ is an irreducible left
topological $A$-module, then it is essential.

Let $A\sb +$ denote the unitization of $A$. Recall that there is
the so-called \emph{canonical morphism} $\pi\sb +\colon A\sb
+\Ptens X\to X$ (resp., $\pi\sb +\colon A\sb +\Ptens X\Ptens A\sb
+\to X$) associated with any left topological $A$-module (resp.,
topological $A$-bimodule) $X$; this morphism is defined by $\pi\sb
+(a\Tens x)=a\cdot x$ (resp., $\pi\sb +(a\Tens x\Tens b)=a\cdot
x\cdot b)$, where $a,b \in A\sb +$, $x\in X$.

Now let us recall some important definitions from the homology theory of
topological algebras (see \cite[Chapter IV]{Hel89},
\cite{Hel93}).

\begin{definition}
A left topological $A$-module (respectively, topological
$A$-bimodule) is \emph{projective} if the canonical
morphism $\pi\sb +$ has a right inverse in the corresponding
category.
\end{definition}

\begin{definition}
A topological algebra $A$ is \emph{left projective} if
the left topological $A$-module $A$ is projective, and
\emph{biprojective} if the topological $A$-bimodule $A$ is
projective.
\end{definition}

\emph{Right projective} topological algebras are defined
similarly. Recall that every biprojective topological algebra is
left and right projective (see \cite[Proposition IV.1.3]{Hel89}).

Recall an important characterization of biprojectivity
for topological algebras.

\begin{proposition}[{see \cite{Hel89,Hel2000}}]\label{proposition3.7}
A topological algebra $A$ is biprojective if and only if the
product map $\pi\sb A\colon A\Ptens A\to A$ has a right inverse in
the category of topological $A$\mbox{-}bimodules.
\end{proposition}

In particular, every biprojective topological algebra $A$ is
idempotent.

\begin{definition}\label{definition3.8}
A topological algebra $A$ is \emph{contractible} if $A$
is biprojective and has an identity.
\end{definition}

The simplest example of a contractible topological algebra is the
$C\sp{\ast}$-algebra $M\sb n(\CC)$ of all complex $n\times
n$-matrices (the full matrix algebra). J.~L.~Taylor noticed in
1972 (see~\cite[p.~181]{Taylor1}; for a proof see \cite[Lemma
11]{Sel96A}) that a topological Cartesian product of full matrix
algebras is a contractible algebra.

Given a left topological $A$-module $X$, the \emph{reduced module}
associated with $X$ is defined by $X\sb {\Pi}=A\ptens{A} X$ (see
\cite{Hel89}). Here $\ptens{A}$ denotes the (complete) projective
tensor product of topological modules over $A$ (see \cite[Chapter
II, \S4.1]{Hel89}).

The following proposition is a generalization of the corresponding
result on Banach modules over Banach algebras (see
\cite[Proposition II.3.13]{Hel89}). For Fr\'echet modules over
Fr\'echet algebras, Proposition \ref{proposition3.9} was
essentially also proved in \cite{Pod2007}.

\begin{proposition}\label{proposition3.9}
Let $A$ be a topological algebra with a left bounded approximate
identity, and let $X$ be a left topological $A$-module. Then the
map
\[
\kappa\sb X\colon X\sb {\Pi}=A\ptens{A} X\to X,\quad a\Tens\sb A
x\mapsto a\cdot x,
\]
is a topological isomorphism onto $\overline{A\cdot X}$.
\end{proposition}

\begin{proof}
Let $L\subset A\ptens{A} X$ denote the linear span of all elements
of the form $a\Tens\limits\sb A x$ ($a\in A$, $x\in X$), and let
$\{e\sb {\nu},\,\nu\in \Lambda\}$  be a left bounded approximate
identity in $A$. If $u=a\Tens\limits\sb A x$, then we have
\[
u=\lim (e\sb \nu a\Tens\sb A x)=\lim (e\sb \nu\Tens\sb A a\cdot
x)=\lim (e\sb \nu\Tens\sb A\kappa\sb X(u)).
\]
Therefore
\begin{equation}
\label{lim-sigma} u=\lim (e\sb \nu\Tens\sb A\kappa\sb
X(u))\quad\text{for each $u\in L$}.
\end{equation}
Given continuous seminorms $p\sb \alpha$ on $A$ and $q\sb \beta$
on $X$ respectively, let $r\sb {\alpha,\beta}$ denote the
corresponding projective tensor seminorm on $A\ptens{A} X$. Then
\eqref{lim-sigma} implies that
\[
r\sb {\alpha,\beta}(u)\le \bigl(\sup p\sb \alpha(e\sb \nu)\bigr)\,
q\sb \beta(\kappa\sb X(u))\quad\text{for each $u\in L$}.
\]
By continuity, the same estimate holds for every $u\in A\ptens{A}
X$. Hence $\kappa\sb X$ is topologically injective%
\footnote[7]{\thinspace In other words, it is a homeomorphism onto
$\Im\kappa\sb X$.}%
. Since the image of $\kappa\sb X$ clearly contains $A\cdot X$ and
is contained in $\overline{A\cdot X}$, the result follows.
\end{proof}

For the next result see \cite[Proposition IV.5.2]{Hel89}.

\begin{proposition}\label{proposition3.10}
Let $A$ be a biprojective topological algebra, and let $X$ be a
left topological $A$-module. Then the reduced module $X\sb {\Pi}$
is projective.
\end{proposition}

From this and from Proposition \ref{proposition3.9} we obtain the
next corollary. Earlier, the corresponding result concerning
Banach algebras and Banach modules was proved in \cite[Theorem
2]{Hel72} (see also \cite[Proposition IV.5.3]{Hel89}).

\begin{corollary}\label{corollary3.11}
Let $A$ be a biprojective topological algebra with a left bounded
approximate identity. Then, for any left topological $A$-module
$X$, the submodule $\overline{A\cdot X}$ is projective.
\end{corollary}

In what follows, given a locally convex space $E$, we denote by
$E\sptilde$ the completion of~$E$. For the reader's convenience,
we recall a result which is a topological version of \cite[Theorem
II.3.17]{Hel89}.

\begin{proposition}[{\cite[Proposition 3.1]{Pir2009}}]
\label{proposition3.12} Let $A$ be a topological algebra, $J$ a
closed left ideal of $A\sb +$, and $X$ a right topological
$A$-module. Then there is a topological isomorphism
\begin{equation}
\label{tens_quot} \alpha\colon X\ptens{A} {}(A\sb
+/J)\sptilde\longrightarrow (X/\overline{X\cdot J})\sptilde
\end{equation}
uniquely determined by $x\Tens\limits\sb A {}(a+J)\mapsto x\cdot
a+\overline{X\cdot J}$.
\end{proposition}

The following result is a generalization of \cite[Lemma
1.1]{Sel79}.

\begin{proposition}
\label{proposition3.13} Let $A$ be a biprojective topological
algebra, and let $J$ be a closed left ideal of $A\sb +$. Then the
left topological $A$-module $(A/\overline{A\cdot J})\sptilde$ is
projective.
\end{proposition}

\begin{proof}
By Proposition \ref{proposition3.12}, the left topological
$A$-module $(A/\overline{A\cdot J})\sptilde$ is topologically
isomorphic to $A\ptens{A} {}(A\sb +/J)\sptilde$. Now the result
follows from Proposition \ref{proposition3.10}.
\end{proof}

\subsection{Homologically trivial and annihilator locally $C\sp{\ast}$-algebras}
The following class of left topological modules will play a
central role in studying homologically trivial locally
$C\sp{\ast}$\mbox{-}algebras.

\begin{definition}[{cf.~\cite{Rief1,Rief2}}]\label{definition3.14}
Let $A$ be a locally $C\sp{\ast}$-algebra. A \emph{Hermitian $A$\mbox{-}module}%
\footnote[8]{\thinspace Such modules are called non-degenerate
star modules in \cite{Hel93} and non-degenerate Hilbert modules in
\cite{Hel98}. We prefer to call them ``Hermitian modules'', as the
phrase ``Hilbert module'' is used in the literature in various
senses (see, e.g.,
\cite{Douglas-Paulsen,Lance,Man-Tro,Phillips1,Gulin}).} is an
essential left topological $A$-module $H$ whose underlying
topological vector space is a Hilbert space, and, moreover,
$\langle a\cdot x,y\rangle=\langle x, a\sp{\ast}\cdot y\rangle$
for all $a\in A$, $x,y\in H$.
\end{definition}

Evidently, the above definition means exactly that the continuous
representation of $A$ associated with our module is a
non-degenerate ${}\sp{\ast}$-representation.

Let us recall the following result.

\begin{theorem}[{\cite[Theorem 5]{Sel77}; see also
\cite[Theorem 4.40]{Sel2000}}] Let $A$ be a $C\sp{\ast}$-algebra.
Then the following conditions are equivalent:
\begin{itemize}
\item[\rm(i)] all irreducible left Banach $A$-modules (or, equivalently%
\footnote[9]{\thinspace As is known (see, e.g.,  \cite[2.9.6 and
2.12.18]{Dixm}), every irreducible left Banach module over a
$C\sp{\ast}$\mbox{-}algebra $A$ is topologically isomorphic
to a Hermitian $A$-module.}%
, all irreducible Hermitian $A$-modules) are projective;
\item[\rm(ii)] $A$ is an annihilator algebra;
\item[\rm(iii)] for every closed left ideal $J$ of $A$, there is a closed left
ideal $L$ of $A$ such that $A=J\oplus L$.
\end{itemize}
\end{theorem}

Following \cite{Barnes1}, we say that an algebra $A$ is a
\emph{modular annihilator algebra} if, for every modular maximal
left ideal $M$ and every modular maximal right ideal $N$ of $A$,
\begin{enumerate}
\item[\rm(i)]
$\ran(M)\neq 0$ and $\ran(A)= 0$, and

\item[\rm(ii)]
 $\lan(N)\neq 0$ and $\lan(A)= 0$.

\end{enumerate}
An equivalent formulation when $A$ is semiprime is as follows: $A$ is a
modular annihilator algebra if and only if $A/\mathrm{Soc}(A)$ is
a radical algebra (see \cite[Theorem 3.4]{Yood}).

Recall from \cite{Barnes3} that a significant number of important
algebras are modular annihilator algebras. In particular, all
semiprime Banach algebras with dense socle, all so-called compact
Banach algebras (see \cite{Alex}) and some algebras of linear
operators are modular annihilator algebras. A well-known result of
B.~A.~Barnes \cite[Theorem 4.2]{Barnes2} (see also \cite[Theorem
8.6.4]{Pal1}) asserts that a semisimple Banach algebra $A$ is
modular annihilator if and only if the spectrum of each element
$a\in A$ has no non-zero accumulation points.

Recall also \cite[Example 4.3]{Yood} that there exists a primitive
(and even topologically simple) modular annihilator Banach algebra
which is not an annihilator algebra. However every modular
annihilator $C\sp{\ast}$\mbox{-}algebra is dual \cite[Theorem
4.1]{Yood}.

For details and references on modular annihilator algebras, see
\cite{Pal1,Pal2}.

\medskip

Now we recall that a Banach space $E$ is said to have the
\emph{approximation property}\footnote[10]{\thinspace This
property is discussed in \cite{Groth,Lind_Tzafr,Hel89,Pal1}. As is
widely known, it was P. Enflo \cite{Enflo} who gave the first
example of a Banach space without the approximation property.} if
the identity operator $\id\sb E$ can be uniformly approximated on
every compact subset $K$ by continuous finite-rank operators
(i.e., for every $\varepsilon>0$ there is a continuous finite-rank
operator $T\colon E\to E$ (depending on $K$ and $\varepsilon$)
satisfying $\|T(x)-x\|\le \varepsilon$ for all $x\in K$).

Recall also the following result.

\begin{theorem}[{\cite[Theorem 4.39]{Sel2000}}]
Let $A$ be a semiprime Banach algebra which satisfies at least one
of the following conditions:
\begin{itemize}
\item[\rm(a)] $A$ has the approximation property;
\item[\rm(b)] all irreducible left Banach $A$-modules have
the approximation property.
\end{itemize}
Then all irreducible left Banach $A$-modules are projective if and
only if $A$ is a modular annihilator algebra.
\end{theorem}

The next result follows immediately from \cite[I, Lemma 19]{Groth}
or \cite[43.2(12)]{Kothe2}.

\begin{lemma}\label{lemma3.17}
Let $E$ be a complete Hausdorff locally convex space, and let $H$
be a Hilbert space. Then for each $u\in E\Ptens H$, $u\neq 0$,
there exist continuous linear functionals $f\in E\sp{\ast}$ and
$g\in H\sp{\ast}$ such that $(f\Ptens g)(u)\neq 0$.
\end{lemma}

The proof of the next proposition, which uses Lemma
\ref{lemma3.17}, is analogous to the proof of
\cite[Lemma~1.4]{Sel79}.

\begin{proposition}\label{proposition3.18}
Let $A$ be a topological algebra, and let $H$ be a projective
Hermitian $A$-module. Then for each $x\in H$, $x\neq 0$, there
exists an $A$-module morphism $\chi\colon H\to A$ such that
$\chi(x)\neq 0$.
\end{proposition}

The next result is \cite[Lemma 4.10.1]{Rick}.

\begin{lemma}\label{lemma3.19}
Let $A$ be an arbitrary ${}\sp{\ast}$-algebra in which $x\sp{\ast}
x=0$ implies $x=0$. Then every minimal left ideal of $A$ is of the
form $Ae$, where $e$ is a unique self-adjoint idempotent. A
similar result holds for right ideals.
\end{lemma}

The next proposition is related to \cite[Theorem 4.10.3]{Rick}.

\begin{proposition}\label{proposition3.20}
Let $A$ be a locally $C\sp{\ast}$-algebra, and let $J$  be a
minimal left ideal of $A$. Then an inner product $\langle
x,\,y\rangle$ can be introduced into $J$ with the following
properties:
\begin{itemize}
\item[\rm(i)] if $p$ is a continuous
$C\sp{\ast}$-seminorm on $A$, then, for each $x\in J$, we have
${p(x)=\lambda\|x\|\sb 0}$, where either $\lambda=0$ or
$\lambda=1$, and $\|x\|\sb 0=\langle x,\,x\rangle\sp {1/2}$, and
so $J$ becomes a Hilbert space with the norm $\|x\|\sb 0$;
\item[\rm(ii)] $\langle ax,\,y\rangle=\langle x,\,a\sp{\ast} y\rangle$ for all $a\in A$, $x,y\in
J$;
\item[\rm(iii)] the left regular representation, $a\mapsto T\sb a$, of $A$
on $J$ is a continuous ${}\sp{\ast}$-representation relative to
$\langle x,\,y\rangle$;
\item[\rm(iv)] the Hermitian $A$-module associated with the left regular
representation of $A$ on~$J$ is projective.
\end{itemize}
\end{proposition}

\begin{proof} By Lemma \ref{lemma3.19}, there exists a self-adjoint
idempotent $e$ such that $J=Ae$. Also, by \cite[Lemma 30.2]{Bon},
$e$ is a minimal idempotent. Since $eAe$ is an Arens-Michael
division algebra, it follows from \cite[Theorem V.1.7]{Hel93} that
$eAe=\mathbb Ce$, i.e., $eAe$ consists of scalar multiples of $e$.
Now, if $x$ and $y$ are any two elements of $J$, then $y\sp{\ast}
x\in eAe$ and hence there exists a scalar $\langle x,\,y\rangle$
such that $y\sp{\ast} x=\langle x,\,y\rangle e$. We shall prove
that $\langle x,\,y\rangle$ is the desired inner product.

It is evident that $\langle x,\,y\rangle$ is linear in
the first variable. Also,
\[\langle y,\,x\rangle e=x\sp{\ast} y=(y\sp{\ast} x)\sp{\ast}=(\langle x,\,y\rangle e)\sp{\ast}
=\overline{\langle x,\,y\rangle} e.
\]
Hence $\langle y,\,x\rangle=\overline{\langle x,\,y\rangle}$. In
particular, we have $\langle x,\,x\rangle=\overline{\langle
x,\,x\rangle}$, i.e., $\langle x,\,x\rangle$ is a real number.
Since
\begin{equation} \label{xxee} x\sp{\ast} x=\langle x,\,x\rangle e\sp{\ast} e,
\end{equation}
it follows from \cite[Corollary 2.5]{Inoue} that $\langle
x,\,x\rangle\ge 0$. Moreover, $\langle x,\,x\rangle=0$ implies
$x\sp{\ast} x=0$, i.e., $x=0$. So $\langle x,\,y\rangle$ is indeed
an inner product.

Now, if $p$ is a continuous $C\sp{\ast}$-seminorm on $A$, then, in
view of (\ref{xxee}), we have
\[p(x)\sp 2=p(x\sp{\ast} x)=\langle x,\,x\rangle p(e\sp{\ast} e)=(\|x\|\sb 0)
\sp 2p(e)\sp 2,\quad x\in J.\] Thus, for each $x\in J$,
$p(x)=\lambda\|x\|\sb 0$, where $\lambda=p(e)$. Since \[p(e)\sp
2=p(e\sp{\ast} e)=p(e),\] it follows that either $\lambda=0$ or
$\lambda=1$. So the topology on the closed subspace $J\subset A$
is generated by the inner product norm $\|x\|\sb 0$. Since $J$ is
complete, $J$~becomes a Hilbert space.

Next, for $a\in A$ and $x,y\in J$, we have \[ \langle
ax,\,y\rangle e=y\sp{\ast} ax=(a\sp{\ast} y)\sp{\ast} x=\langle
x,\,a\sp{\ast} y\rangle e,
\]
and hence
\begin{equation} \label{ax-ye}
\langle ax,\,y\rangle=\langle x,\,a\sp{\ast} y\rangle.
\end{equation}

Moreover, in view of (\ref{ax-ye}), we have \[\langle T\sb
ax,\,y\rangle=\langle x,\,T\sb {a\sp{\ast}}y\rangle,\] and thus we
obtain $T\sb {a\sp{\ast}}=T\sb a\sp{\ast}$, which shows that
$a\mapsto T\sb a$ is a ${}\sp{\ast}$-representation on $J$.

Finally, let $H$ be the Hermitian $A$-module associated with the
left regular representation of $A$ on $J$. Consider the operator
\[\rho\colon H\to A\sb +\Ptens H,\quad x\mapsto x\Ptens e,\]
where $x=ae\in H=J\subset A\sb +$. It is obvious that $\rho$ is an
$A$-module morphism such that $\pi\sb +\circ \rho$ is the identity
on $H$, where $\pi\sb +\colon A\sb +\Ptens H\to H$ is the
canonical morphism. It follows that $H$
is projective.
\end{proof}

The next result characterizes projective irreducible Hermitian
modules over a locally $C\sp{\ast}$-algebra.

\begin{corollary}\label{corollary3.21}
Let $H$ be an irreducible Hermitian module over a locally
$C\sp{\ast}$-algebra. Then the following conditions are
equivalent:
\begin{itemize}
\item[\rm(i)] $H$ is projective;
\item[\rm(ii)] $H$ is topologically isomorphic to a minimal left ideal of $A$.
\end{itemize}
\end{corollary}

\begin{proof} (i)$~\Rightarrow~$(ii) Since $H$ is projective, it follows from
Proposition \ref{proposition3.18} that there exists an $A$-module
morphism $\chi\colon H\to A$ such that $\chi(x)\neq 0$ for some
$x\neq 0$. Then $\Ker\chi\neq H$ and $\Im\chi\neq 0$. Since $H$ is
irreducible, it follows at once that $\Ker\chi=0$ and $\Im\chi$ is
an irreducible submodule of $A$. The latter implies that $\Im\chi$
is a minimal left ideal of~$A$. Since, by Proposition
\ref{proposition3.20}, every minimal left ideal of $A$ is a
Hilbert space, the rest follows from the Open Mapping Theorem.

(ii)$~\Rightarrow~$(i) This follows from Proposition
\ref{proposition3.20}.
\end{proof}

The following result is a part of \cite[Proposition 8.4.4]{Pal1}.

\begin{lemma}\label{lemma3.22}
Let $A$ be a semiprime algebra and let $e$ be a minimal idempotent
of $A$. Then $A(1-e)$ is a modular maximal left ideal of $A$.
\end{lemma}

\begin{definition}\label{definition3.23}
Let $A$ be a topological algebra. We call $A$ a \emph{topologically
modular annihilator algebra} if, for every closed modular maximal
left ideal $M$ and every closed modular maximal right ideal $N$ of
$A$,
\begin{enumerate}
\item[\rm(i)]
$\ran(M)\neq 0$ and $\ran(A)= 0$, and

\item[\rm(ii)]
$\lan(N)\neq 0$ and $\lan(A)= 0$.
\end{enumerate}
\end{definition}

We recall (cf. \cite[Chapter VI, \S1.2]{Hel93}) that a closed
submodule $X\sb 0$ of a left (right, or two-sided) topological
module $X$ over a topological algebra $A$ is said to be
\emph{complemented as a topological module} if there exists
another closed submodule $X\sb 1$ of $X$ such that the operator
\[\lambda\colon X\sb 0\times X\sb 1\to X,\quad (x,\,y)\mapsto x+y,\]
is an isomorphism of topological modules.

Now we are in a position to answer the first question posed in the
beginning of the paper.

\begin{theorem}\label{theorem3.24}
Let $A$ be a locally $C\sp{\ast}$-algebra. Then the following
conditions are equivalent:
\begin{itemize}
\item[\rm(i)] all irreducible Hermitian $A$-modules are projective;
\item[\rm(ii)] $\mathrm{Soc}(A)$ is dense in $A$;
\item[\rm(iii)] $A$ is the direct topological sum of its minimal closed two-sided
ideals each of which is topologically ${}\sp{\ast}$-isomorphic to
an elementary $C\sp{\ast}$-algebra;
\item[\rm(iv)] $A$ is an annihilator algebra;
\item[\rm(v)] for every closed left ideal $J$ of
$A$ and for every closed right ideal $K$ of $A$, we have
$A=J\oplus {}(\ran(J))\sp{\ast}$ and $A=K\oplus {}
(\lan(K))\sp{\ast}$, and moreover $J$ and $K$ are complemented as
topological submodules of $A$;
\item[\rm(vi)] $A$ is a dual algebra;
\item[\rm(vii)] $A$ is a complemented algebra;
\item[\rm(viii)] $A$ is a left quasi-complemented algebra;
\item[\rm(ix)] for every closed modular maximal left ideal $M$ of
$A$, there is a closed left ideal $L$ of~$A$ such that $A=M\oplus
L$;
\item[\rm(x)] the closed modular maximal left ideals of
$A$ are precisely the left ideals of the form $A(1-e)$,
where $e$ is a minimal idempotent%
\footnote[11]{\thinspace Obviously, such left ideals can also be
characterized as those of the form $\lan (R)$ for some minimal
right ideal $R$ of $A$.} of $A$;
\item[\rm(xi)] all irreducible barrelled%
\footnote[12]{\thinspace Recall \cite{Bourbaki} that a locally
convex space $E$ is {\em barrelled} if every closed absorbing
absolutely convex subset of $E$ is a neighbourhood of zero. In
particular, all Fr\'echet spaces are barrelled.} left topological
$A$-modules are projective;
\item[\rm(xii)] $A$ is a topologically modular annihilator algebra.
\end{itemize}
\end{theorem}

\begin{proof} (i)$~\Rightarrow~$(ii) Let
$I=\overline{\mathrm{Soc}(A)}\neq A$. Then $I$ is a proper closed
two-sided ideal of~$A$. By \cite[Theorem 11.7]{Frag3}, $I$ is a
${}\sp{\ast}$-ideal and the quotient algebra $A/I$ equipped with
the
quotient topology is a $C\sp{\ast}$-convex algebra%
\footnote[13]{\thinspace In other words, a (not necessarily
complete) locally $m$-convex ${}\sp{\ast}$-algebra whose topology
is defined
by a family of $C\sp{\ast}$-seminorms.}%
. According to \cite[Corollary 20.4]{Frag3}, every
$C\sp{\ast}$\mbox{-}convex algebra has enough continuous
topologically irreducible ${}\sp{\ast}$-representations to
separate its points. Therefore there exists a non-zero continuous
topologically irreducible ${}\sp{\ast}$-representation $T\colon
A/I\rightarrow \B(H)$ on some Hilbert space $H$. (Here, as usual,
$\B(H)$ is the Banach algebra of all continuous linear operators
on $H$.) Consider the quotient map $\sigma\colon A\rightarrow A/I$
and the topologically irreducible
${}\sp{\ast}$\mbox{-}representation $\widetilde{T}=T\circ \sigma$
of $A$ on $H$. Then, by \cite[Theorem 19.2]{Frag3},
$\widetilde{T}\colon A\rightarrow \B(H)$ is algebraically
irreducible. The associated Hermitian $A$-module $H$ is
irreducible and hence projective. By Proposition
\ref{proposition3.18}, there exists an $A$-module morphism
$\chi\colon H\to A$ such that $\chi(x)\neq 0$ for some $x\neq 0$.
Then $\Ker\chi\neq H$ and $\Im\chi\neq 0$. Let $J=\Im\chi$. Since
$H$ is irreducible, it follows that $\Ker\chi=0$ and $J$ is a
minimal left ideal of $A$. Hence $J\subset \mathrm{Soc}(A)\subset
I$. Since, obviously, $\widetilde{T}(I)=0$, we have $I\cdot H=0$
and \[J\sp 2\subset I\cdot J\subset \chi(I\cdot H)=0.\] Since the
algebra $A$ is semiprime, we get $J=0$. So we have a
contradiction.

(ii)$~\Rightarrow~$(iii) This follows from Theorem
\ref{theorem3.3}.

(iii)$~\Rightarrow~$(iv) This is essentially Theorem
\ref{theorem3.3}(v).

(iv)$~\Rightarrow~$(v) This follows from Lemma \ref{lemma2.3} and
from inequalities (\ref{pcltpa}) and (\ref{pbltpa}) (see the
proof of Lemma \ref{lemma2.3}).

(v)$~\Rightarrow~$(vi) Let $J$ be a closed left ideal of $A$. By
condition (v), we have
\begin{equation} \label{aranjoplusj}
A=J\oplus {}(\ran(J))\sp{\ast}.
\end{equation}
Since $\ran(J)$ is a closed right ideal of $A$, it follows from
condition (v) that \[A=\ran(J)\oplus {}(\lan(\ran(J)))\sp{\ast}.\]
Applying the involution to this equality, we obtain
$A=(\ran(J))\sp{\ast}\oplus \lan(\ran(J))$, i.e.,
\begin{equation} \label{aranjopluslrj}
A=\lan(\ran(J))\oplus {}(\ran(J))\sp{\ast}.
\end{equation}
Since $J\subset\lan(\ran(J))$, it follows from (\ref{aranjoplusj})
and (\ref{aranjopluslrj}) that $\lan(\ran(J))=J$. The proof of the
equality ${\ran(\lan(K))=K}$, for every  closed right ideal $K$ of
$A$, is similar.

(vi)$~\Rightarrow~$(vii) For every closed left ideal $J$ of $A$,
we set $J\sp q=(\ran(J))\sp{\ast}$. Then it is clear that $J\sp q$
is a closed left ideal of $A$ and $J\sp q=\lan(J\sp{\ast})$.
Moreover, we have
\[(J\sp q)\sp q=\lan((J\sp q)\sp{\ast})=\lan((\ran(J))\sp {\ast\ast})=\lan(\ran(J))=J,\]
because $A$ is dual by condition (vi).

Now note that $A$, being a dual algebra, is an annihilator
algebra. From this and from Lemma \ref{lemma2.3} we get that
$A=J\oplus J\sp q$ for every closed left ideal $J$.

Finally, it is evident that, if $J\sb 1$ and $J\sb 2$ are closed
left ideals and if $J\sb 1\subset J\sb 2$, then $J\sb 2\sp
q\subset J\sb 1\sp q$. Thus $A$ is a left complemented algebra
with left complementor $q\colon J\mapsto J\sp q$.

Similarly, it can be proved that $A$ is a right complemented
algebra with right complementor $p\colon K\mapsto K\sp p$, where
$K\sp p=(\lan(K))\sp{\ast}$ for every closed right ideal $K$.

(vii)$~\Rightarrow~$(viii) This is trivial.

(viii)$~\Rightarrow~$(ix) Let $q\colon J\mapsto J\sp q$ be a left
quasi-complementor on $A$. Then, since \[A\sp q=A\sp q\cap A=0,\]
we have $0\sp q=(A\sp q)\sp q=A$. Now let $M$ be an arbitrary
closed modular maximal left ideal of $A$. Then $M\sp q\neq 0$,
since otherwise
\[M=(M\sp q)\sp q=0\sp q=A.\]
Since $M+M\sp q$ is a left ideal which contains $M$ properly, it
follows that $M+M\sp q=A$. Putting $L=M\sp q$, we get $A=M\oplus
L$.

(ix)$~\Rightarrow~$(x) By Lemma \ref{lemma3.22}, we just need to
see that any closed modular maximal left ideal $M$ of $A$ is of
the form  $A(1-e)$ for a minimal idempotent $e\in A$. By condition
(ix), for every such left ideal $M$ there is a closed left ideal
$L$ of $A$ such that $A=M\oplus L$. Since the left ideal $M$ is
maximal, it follows that $L$ is minimal.

Let $u$ be a right identity for $A$ modulo $M$. Since $A=M\oplus
L$, we can write $u=m+e$, where $m\in M$ and $e\in L$. It is clear
that $A(1-e)\subset M$, i.e., $e$ is also a right identity for $A$
modulo $M$. Let now $b\in L$. Since
\[b-be\in M\cap L=0,\] we have $b=be$. It follows that $L=Ae$ and
$e\sp 2=e$. Since $L$ is a minimal left ideal of $A$, the
idempotent $e$ is minimal (see \cite[Lemma 30.2]{Bon}).

Suppose now that $M$ properly includes $A(1-e)=\{a\in A\mid
ae=0\}$. Then there is some $a\in M$ with $ae\neq 0$. Since
$a-ae\in M$, we have \[ae\in M\cap Ae=M\cap L=0.\] This
contradiction shows that $M=A(1-e)$.

(x)$~\Rightarrow~$(xi) Let $X$ be an irreducible barrelled left
topological $A$-module. Taking an arbitrary element ${x\in
X\setminus \{0\}}$ and using \cite[Theorem VI.2.8]{Hel93} we get
that the set
\[M=\{a\in A \mid a\cdot x=0\}\] is a closed modular maximal left
ideal of $A$. By condition (x), we have $M=A(1-e)$, where $e\in A$
is a minimal idempotent. Since, for $a\in A$, we have $a-ae\in M$
and hence
\[a\cdot (x-e\cdot x)=(a-ae)\cdot x=0,\] we get that $e\cdot x=x$
(see \cite[Proposition VI.2.6(II)]{Hel93}).

Let $L=Ae$. By \cite[Proposition~30.6]{Bon}, $L$ is a minimal left
ideal of $A$. It follows from Proposition \ref{proposition3.20}
that $L$ is a Hilbert space. Furthermore, it is a projective
irreducible Hermitian $A$-module. Consider the operator
\[\varphi\colon L\to X,\quad a\mapsto a\cdot x.\]
Clearly, $\varphi$ is an $A$-module morphism such that
\[\Ker\varphi=M\cap L=0\] and $\varphi(e)=x\neq 0$.
Since $X$ is irreducible, it follows that $\Im\varphi=X$. Since
$X$ is barrelled, the Open Mapping Theorem \cite[Chapter
VI, Theorem 7 and Corollary~2]{rob-eng} implies that $\varphi$ is
a topological isomorphism. Thus $X$ is topologically isomorphic to
a projective left topological $A$-module, and so it is projective.

(xi)$~\Rightarrow~$(i) This is trivial, since all Hilbert spaces
are barrelled.

(x)$~\Rightarrow~$(xii) Let $M$ be a closed modular maximal left
ideal of $A$. Since, by condition~(x), we have $M=A(1-e)$, where
$e$ is a minimal idempotent of $A$, it follows that $e\in
\ran(M)$. Hence $\ran(M)\neq 0$. Using the involution it is easily
seen that every closed modular maximal right ideal $N$ of $A$ has
a non-zero left annihilator. Thus $A$ is a topologically modular
annihilator algebra.

(xii)$~\Rightarrow~$(ix) Let $M$ be a closed modular maximal left
ideal of $A$. By condition (xii), we have $\ran(M)\neq 0$. Define
\[L=(\ran(M))\sp{\ast}.\] Then $L$ is a non-zero closed left ideal
of~$A$. By Lemma \ref{lemma2.3}(i), we have $M\cap L=0$. Since
$M\oplus L$ is a left ideal which contains $M$ properly, it
follows that $M\oplus L=A$.
\end{proof}

\begin{remark}
The implications (iv)~$\Rightarrow$~(vi)~$\Rightarrow$~(vii) in
Theorem \ref{theorem3.24} were proved earlier in \cite{Har2005}.
We included a proof for completeness.
\end{remark}

\begin{remark}
Within the setting of $C\sp{\ast}$-algebras, many other conditions
equivalent to those listed in Theorem \ref{theorem3.24} were
obtained in \cite{Frank1,Frank2}.
\end{remark}

From Theorems \ref{theorem3.3} and \ref{theorem3.24} we
immediately get the following result.

\begin{proposition}\label{proposition3.26}
Let $A$ be an annihilator locally $C\sp{\ast}$-algebra, and let
$I$ be a closed two-sided ideal of $A$. Then $I$ and $A/I$ are
annihilator locally $C\sp{\ast}$-algebras, and moreover $A$ is
topologically ${}\sp{\ast}$-isomorphic to the Cartesian product of
the algebras $I$ and $A/I$.
\end{proposition}

Propositions \ref{proposition3.1} and \ref{proposition3.26} imply
the following corollary.

\begin{corollary}
Let $A$ be an annihilator locally $C\sp{\ast}$-algebra, and let
$P$ be the family of all continuous $C\sp{\ast}$-seminorms on $A$.
Then, for each $p\in P$, the concomitant $C\sp{\ast}$-algebra
$A\sb p$ is annihilator and, as a consequence, $A$ can be
represented as an inverse limit of annihilator
$C\sp{\ast}$\mbox{-}algebras.
\end{corollary}

Lemma \ref{lemma2.4} and Proposition \ref{proposition3.26} imply
the following.

\begin{proposition}\label{proposition3.28}
The Cartesian product $A=\prod\limits\sb {\nu\in \Lambda} A\sb
{\nu}$ of a family $\{A\sb {\nu}\mid \nu\in \Lambda\}$ of locally
$C\sp{\ast}$-algebras, with the product topology, is an
annihilator algebra if and only if all the algebras $A\sb {\nu}$
are annihilator.
\end{proposition}

\subsection{Unital, metrizable, and non-unital annihilator locally $C\sp{\ast}$-algebras}
Theorem \ref{theorem3.24} can be strengthened in certain cases.
First note that Corollary \ref{corollary3.4}, Theorem~\ref{theorem3.24} and Proposition~\ref{proposition3.28} imply the
following result.

\begin{theorem}\label{theorem3.29}
Let $A$ be a unital locally $C\sp{\ast}$-algebra. Then $A$
satisfies the equivalent conditions of Theorem \ref{theorem3.24}
if and only if $A$ is topologically ${}\sp{\ast}$-isomorphic to
the Cartesian product of a family of full matrix algebras.
\end{theorem}

Next we recall that a $\sigma$-$C\sp{\ast}$-algebra is a
metrizable locally $C\sp{\ast}$-algebra (or, equivalently, a
locally $C\sp{\ast}$-algebra whose topology is determined by a
countable family of $C\sp{\ast}$-seminorms). Every
$\sigma$-$C\sp{\ast}$-algebra is topologically
${}\sp{\ast}$-isomorphic to the inverse limit of a sequence
\begin{equation}
\label{proj1} A\sb 1 \xleftarrow{\sigma\sb 1} A\sb 2
\xleftarrow{\sigma\sb 2} A\sb 3\xleftarrow{\sigma\sb 3}\cdots
\end{equation}
of $C\sp{\ast}$-algebras and surjective
${}\sp{\ast}$-homomorphisms $\sigma\sb n\colon A\sb {n+1}\to A\sb
n$ (see \cite[Section 5]{Phillips1}). Moreover, all the natural
homomorphisms $A\to A\sb n\; (n\in \mathbb{N})$ are also
surjective.

\begin{theorem}\label{theorem3.30}
Let $A$ be an annihilator $\sigma$-$C\sp{\ast}$-algebra. Then $A$
is topologically ${}\sp{\ast}$-iso\-mor\-phic to the Cartesian
product of a countable family of annihilator
$C\sp{\ast}$-algebras.
\end{theorem}

\begin{proof}
Choose an inverse system \eqref{proj1} of $C\sp{\ast}$-algebras
such that
\[
A\cong\varprojlim A\sb n.
\]
As was noted above, we may assume that all the maps $A\sb {n+1}\to
A\sb n$ and $A\to A\sb n$ are onto. By the Open Mapping Theorem,
$A\sb n$ is topologically ${}\sp{\ast}$-isomorphic to $A/I\sb n$,
where
\[
I\sb n=\Ker {}(A\to A\sb n).
\]
It follows from Proposition \ref{proposition3.26} that $A\sb n$ is
an annihilator $C\sp{\ast}$-algebra for each $n\in\mathbb{N}$. So,
by \cite[4.7.20]{Dixm}, each $A\sb n$ is isometrically
${}\sp{\ast}$-isomorphic to a $c\sb 0$-sum of elementary
$C\sp{\ast}$\mbox{-}algebras. In particular, the spectrum
$\widehat{A}\sb n$ of $A\sb n$ is discrete. Now \cite[Proposition
5.2]{PiS2007} implies that there exists a family $\{B\sb n\mid
n\in\mathbb{N}\}$ of $C\sp{\ast}$-algebras such that $A$ is
topologically ${}\sp{\ast}$\mbox{-}isomorphic to $\prod\limits\sb
n B\sb n$. Using Proposition \ref{proposition3.28}, we see that
each $B\sb n$ is an annihilator algebra. This completes the proof.
\end{proof}

By combining this result with Theorem \ref{theorem3.24} and
Proposition \ref{proposition3.28} we get the following.

\begin{corollary}\label{corollary3.31}
Let $A$ be a $\sigma$-$C\sp{\ast}$-algebra. Then $A$ satisfies the
equivalent conditions of Theorem \ref{theorem3.24} if and only if
$A$ is topologically ${}\sp{\ast}$-isomorphic to the countable
Cartesian product $\prod\limits\sb n A\sb n$, where each $A\sb n$
is a $C\sp{\ast}$-algebra isomorphic to a $c\sb 0$\mbox{-}sum of
elementary $C\sp{\ast}$-algebras.
\end{corollary}

\begin{remark}
Theorem \ref{theorem3.30} and Corollary~\ref{corollary3.31} can
not be extended to arbitrary (i.e., not necessarily metrizable)
locally $C\sp{\ast}$-algebras. Namely, there exists a non-unital
annihilator locally $C\sp{\ast}$-algebra that is not topologically
isomorphic to a Cartesian product of annihilator
$C\sp{\ast}$-algebras (see Theorem \ref{theorem4.28} below).
\end{remark}

At the same time we have the following result.

\begin{proposition}\label{proposition3.33}
Let $A$ be a non-unital annihilator locally $C\sp{\ast}$-algebra.
Then $A$ is topologically ${}\sp{\ast}$-isomorphic to the
Cartesian product of two annihilator locally $C\sp{\ast}$-algebras
one~of which is an infinite-dimensional $C\sp{\ast}$-algebra.
\end{proposition}

\begin{proof}
By Theorem \ref{theorem3.24}, $\mathrm{Soc}(A)$ is dense in $A$.
Since $A$ is non-unital, it follows from Corollary
\ref{corollary3.4}(ii) that there exists a continuous
$C\sp{\ast}$-seminorm $p$ on $A$ such that the concomitant
$C\sp{\ast}$-algebra $A\sb p$ is infinite-dimensional. We set
$I=\Ker p$. Then $I$ is a closed two-sided ideal of $A$ and, by
Proposition~\ref{proposition3.1}, $A/I$ is topologically
${}\sp{\ast}$-isomorphic to the $C\sp{\ast}$-algebra $A\sb p$. It
remains to apply Proposition \ref{proposition3.26}.
\end{proof}

\section{Biprojective locally $C\sp{\ast}$-algebras}
\subsection{Examples and general properties of biprojective algebras}
We recall that biprojective Banach and topological algebras were
introduced by A.~Ya.~Helemskii \cite{Hel72,Hel73} and have been
studied by many authors (see, e.g.,
\cite{Hel72,Hel72B,Joh72,Sel75,Sel76,Sel77,Hel78,Sel79,Lyk79,Rach81,Sel92,Sel98,Sel99,Sel2000,Sel2002,Aris2002,Pir2002,Wood02,Pir2004,Pir2004a,PiS2007,Aris2008,Runde_Fourier}).

The original motivation to study such algebras was the vanishing
of their cohomology groups%
\footnote[14]{\thinspace The definition of cohomology groups of
topological algebras can be found, e.g., in
\cite{Hel89}.}%
, ${\mathcal{H}\sp n(A,X)}$, with coefficients in arbitrary
topological $A$-bimodules $X$ for all $n\ge 3$ (see, e.g.,
\cite[Theorem 2.4.21]{Hel2000}). The structure of biprojective
semisimple Banach algebras with the approximation property is
described in \cite{Sel79} (see also \cite{Hel89} and
\cite{Runde}). In~\cite{Sel99}, the cohomology groups of
biprojective Banach algebras are completely computed for
arbit\-rary coefficients. The description of these groups is given
in terms of double multipliers and quasi-multipliers of a given
bimodule of coefficients. As an application, biprojective Banach
algebras are characterized in terms of their cohomology groups.
In~particular, it is shown (see \cite[Theorem 5.9]{Sel99}) that a
Banach algebra is biprojective if and only if its one-dimensional
cohomology groups with coefficients in bimodules of double
multipliers are trivial.

We recall some examples of biprojective algebras.

\begin{example}\label{example4.1}
As was noted after Definition \ref{definition3.8}, the Cartesian
product $\prod\limits\sb {\nu\in\Lambda} M\sb {n\sb {\nu}}(\CC)$
of any family of full matrix algebras is a contractible
topological algebra, i.e, a biprojective algebra with an identity.
Moreover, we recall that the Cartesian product of any family of
contractible topological algebras is a contractible algebra (see
\cite[Lemma 11]{Sel96A}).
\end{example}

Here it is important to note that a topological algebra $A$ is
contractible if and only if its cohomology groups with
coefficients in arbitrary topological $A$-bimodules vanish for all
$n\ge 1$ (see \cite[Theorem IV.5.8]{Hel89}). For some results
concerning contractible algebras, see
\cite{Taylor1,Taylor2,Liddell,Sel76A,Sel77,Hel89,Sel96A,Runde1,Frag2,Sel2002,Pir2004,PiS2007}.

\begin{example}
The simplest example of a non-contractible biprojective Banach
algebra is perhaps the Banach sequence algebra $\ell\sb 1$ with
coordinatewise multiplication. The Banach sequence algebra $c\sb
0$ is also biprojective (see \cite[Examples
IV.5.9--IV.5.10]{Hel89} or \cite[Example VII.1.80]{Hel93}).
\end{example}

\begin{example}
If $G$ is a compact group, then the group algebras $L\sp 1(G)$ and
$C\sp{\ast}(G)$ are biprojective Banach algebras (\cite{Hel72},
see also \cite{Hel89}).
\end{example}

\begin{example}\label{example4.4}
The $c\sb 0$-sum, $\bigoplus\nolimits\sb 0\{M\sb {n\sb
{\nu}}(\CC)\mid \nu\in \Lambda\}$, of a family of full matrix
$C\sp{\ast}$\mbox{-}algebras is a biprojective
$C\sp{\ast}$-algebra \cite[Theorem 3]{Hel72}.
\end{example}

Actually, there are no other biprojective $C\sp{\ast}$-algebras:

\begin{theorem}[{\cite{Sel76,Sel77};
see also \cite[Theorem 4.62]{Sel2000}}] \label{theorem4.5} Every
biprojective $C\sp{\ast}$-algebra is isometrically
${}\sp{\ast}$-isomorphic to the $c\sb 0$-sum of a family of full
matrix algebras. In~particular, every commutative biprojective
$C\sp{\ast}$-algebra is isometrically ${}\sp{\ast}$-isomorphic to
a
$C\sp{\ast}$-algebra of the form%
\footnote[15]{\thinspace As usual (see, e.g., \cite{Lind_Tzafr}),
we write $c\sb 0(\Lambda)$ for $\bigoplus\nolimits\sb 0\{A\sb
\nu\mid \nu\in \Lambda\}$, where $A\sb \nu=\CC$ for each $\nu\in
\Lambda$.} $c\sb 0(\Lambda)$ for some set $\Lambda$.
\end{theorem}

\begin{example}
The Cartesian product of a countable family of biprojective
Fr\'echet algebras is a biprojective Fr\'echet algebra
\cite[Proposition 1.15]{PiS2007}. As a consequence, the Cartesian
product of a countable family of biprojective
$C\sp{\ast}$-algebras is a biprojective
$\sigma$\mbox{-}$C\sp{\ast}$\mbox{-}algebra.
\end{example}

Actually, there are no other biprojective
$\sigma$\mbox{-}$C\sp{\ast}$\mbox{-}algebras:

\begin{theorem}[{\cite[Theorem 5.3]{PiS2007}}]\label{theorem4.7}
Every biprojective $\sigma$-$C\sp{\ast}$-algebra is topologically
${}\sp{\ast}$-isomorphic to the countable Cartesian product
$\prod\limits\sb n A\sb n$, where each $A\sb n$ is a
$C\sp{\ast}$-algebra isomorphic to a $c\sb 0$\mbox{-}sum of full
matrix algebras.
\end{theorem}

\begin{example}
Let $\alpha=(\alpha\sb n)$, $n\in\N$, be a non-decreasing sequence
of positive numbers with $\lim\limits\sb {n\to \infty}\alpha\sb
n=+\infty$. Then the power series spaces $\Lambda\sb 1(\alpha)$
and $\Lambda\sb \infty(\alpha)$ are biprojective Fr\'echet
algebras under pointwise product (see \cite{Pir2002} or
\cite{Pir2004}). In particular, the algebra $s$ of rapidly
decreasing sequences is biprojective.
\end{example}

\begin{example}\label{example4.9}
If $G$ is a compact Lie group, then the Fr\'echet algebra
$\mathcal E(G)$ of smooth functions on $G$ with convolution
product is biprojective (\cite[Example 6]{Sel96B}; see also
\cite{Pir2002}).  Note that $\mathcal E(G)$ is not a Banach
algebra unless $G$ is finite.
\end{example}

\begin{example}
If $G$ is a compact Lie group, then the algebra $\mathcal{E}'(G)$
of distributions on $G$ with convolution multiplication is a
biprojective topological algebra with an identity
(\cite[Proposition 7.3]{Taylor2}, see also \cite[Assertion
IV.5.30]{Hel89}). In other words, $\mathcal{E}'(G)$ is a
contractible algebra. Note that the biprojective Fr\'echet algebra
$\mathcal E(G)$ from Example \ref{example4.9} is a dense ideal in
$\mathcal{E}'(G)$.
\end{example}

\begin{example}\label{example4.11}
Let $(E,F)$ be a pair of complete Hausdorff locally convex spaces,
and let
$\langle\:\cdot\: ,\:\cdot\:\rangle\colon E\times F\to\mathbb{C}$
be a jointly continuous bilinear form
that is not identically zero. The space $A=E\Ptens F$ is then a
topological algebra with respect to the multiplication defined by
\[
(x\sb 1\Tens y\sb 1)(x\sb 2\Tens y\sb 2)=\langle x\sb 2,y\sb
1\rangle\, x\sb 1\Tens y\sb 2 \quad (x\sb i\in E,\; y\sb i\in F).
\]
This algebra is biprojective (\cite[Example 4]{Sel96B}, see also
\cite{Sel79}, \cite{Pir2002} and \cite{Pir2004}).

In particular, if $E$ is a Banach space with the approximation
property, then the algebra $A=E\Ptens E\sp{\ast}$ is isomorphic to
the Banach algebra $\mathcal N(E)$ of all nuclear operators on
$E$, and so $\mathcal N(E)$ is biprojective.
\end{example}

Actually, there are no other biprojective Banach algebras of the
form $\mathcal N(E)$:

\begin{theorem}[{see \cite[Corollaries 1, 3]{Sel92}}]
Let $E$ be a Banach space, and let $A=\mathcal N(E)$. Then the
following conditions are equivalent:
\begin{itemize}
\item[\rm(i)] $A$ is a biprojective algebra;
\item[\rm(ii)] ${\mathcal{H}\sp n(A,X)}=0$ for all Banach $A$-bimodules $X$ and for all $n\ge 3$;
\item[\rm(iii)] there exists $k\ge 3$ such that ${\mathcal{H}\sp n(A,X)}=0$
for all Banach $A$-bimodules $X$ and for all $n\ge k$;
\item[\rm(iv)] $E$ has the approximation
property.
\end{itemize}
\end{theorem}

Let us briefly explain why the algebra $\mathcal N(E)$ cannot be
biprojective unless $E$ has the approximation property. It is
known that, for a biprojective Banach algebra $A$, the operator
\[
\kappa\sb A\colon A\ptens{A}A\to A,\quad a\Tens\sb Ab\mapsto ab,
\]
is a topological isomorphism (see \cite[Corollary 4.2]{Sel99}).
Let now $A=\mathcal N(E)$, where $E$ is a Banach space.
By \cite[Lemma 3]{Sel92}, the $A$-bimodule $A\ptens{A}A$ is isomorphic to $E\Ptens E\sp{\ast}$,
and the operator $\kappa\sb A$ can be identified with the
so-called \emph{trace homomorphism}
\[
\mathrm{Tr}\colon E\Ptens E\sp{\ast}\to \mathcal N(E),\quad
\mathrm{Tr}(x\Tens f)(y)=\langle y,f\rangle\, x\qquad (x,y\in E,\:
f\in E\sp{\ast}),
\]
of the dual pair $(E,E\sp{\ast})$ (see \cite{Hel89}). However,
$\mathrm{Tr}$ is known to have a non-trivial kernel unless $E$ has
the approximation property (see \cite[I, Proposition 35]{Groth}).

Let us also note that, for any Banach space $E$, the Banach algebra
$\mathcal N(E)$ has the weaker property of being
``quasi-biprojective''. For details, see \cite{Sel2004}.

\medskip

The following result of O.~Yu.~Aristov describes the structure of
finite-dimensional biprojective algebras.

\begin{theorem}[{see \cite[Theorem 7.1]{Aris2008}}]\label{theorfinite}
Let $A$ be a finite-dimensional biprojective algebra. Then $A$ is
isomorphic to the Cartesian product of finitely many algebras of
the type $E\Ptens F$, where
$(E,F,\langle\:\cdot\:,\:\cdot\:\rangle )$ is a pair of
finite-dimensional spaces together with a non-zero bilinear form.
\end{theorem}

We pass from examples to some general facts. Recall that, if $A$
and $B$ are topological algebras, then $A\Ptens B$ is also a
topological algebra with multiplication defined by
\[
(a\sb 1\Tens b\sb 1)(a\sb 2\Tens b\sb 2)=a\sb 1a\sb 2\Tens b\sb
1b\sb 2 \qquad (a\sb 1,a\sb 2\in A,\; b\sb 1,b\sb 2\in B).
\]
By $A\sp{\rm op}$ we denote the topological algebra
\emph{opposite} to $A$ (with the same underlying space as $A$, but
with multiplication $a\circ b$ equal to the ``previous'' $ba$).
The next proposition follows immediately from Proposition
\ref{proposition3.7}.

\begin{proposition}
Let $A$ and $B$ be topological algebras.
\begin{itemize}
\item[\rm(i)] If $A$ and $B$ are biprojective, then so is $A\Ptens B$.
\item[\rm(ii)] If $A$ is biprojective, then so is
$A\sp {\rm op}$.
\end{itemize}
\end{proposition}

The following is a generalization of the corresponding result on
Banach algebras (see \cite[Lemma 1.3]{Sel79}).

\begin{proposition}
\label{proposition4.14} Let $A$ be a biprojective topological
algebra, and let $I$ be a closed two-sided ideal of $A$. Then
$(A/\overline{A\cdot I})\sptilde$ is a biprojective topological
algebra.
\end{proposition}

\begin{proof}
Set $B=(A/\overline{A\cdot I})\sptilde$, and consider the
commutative diagram
\[
\xymatrix{
A\Ptens A \ar[r]\sp (.6){\pi\sb A} \ar[d]\sb {\sigma\Ptens \sigma} & A \ar[d]\sp \sigma \\
B \Ptens B \ar[r]\sp (.6){\pi\sb B} & B }
\]
in the category of topological $A$\mbox{-}bimodules. Here $\pi\sb
A$ and $\pi\sb B$ are the product maps, and $\sigma$ is the
natural operator from $A$ to $(A/\overline{A\cdot I})\sptilde$.
Using Proposition \ref{proposition3.7}, choose an $A$-bimodule
morphism $\rho\sb A\colon A\to A\Ptens A$ such that $\pi\sb
A\circ\rho\sb A=\id\sb A$. We have
\[
\rho\sb A(A\cdot I)\subset (A\Ptens A)\cdot
I\subset\Ker(\sigma\Ptens \sigma),
\]
and so there exists a unique $A$-bimodule morphism $\rho\sb B$
making the diagram
\[
\xymatrix{
A\Ptens A \ar[d]\sb {\sigma\Ptens \sigma} & A \ar[d]\sp \sigma \ar[l]\sb (.4){\rho\sb A} \\
B \Ptens B & B \ar[l]\sp (.4){\rho\sb B} }
\]
commutative. Therefore,
\[
\pi\sb B\circ\rho\sb B\circ\sigma=\pi\sb B\circ(\sigma\Ptens
\sigma)\circ\rho\sb A= \sigma\circ\pi\sb A\circ\rho\sb A=\sigma.
\]
Since $\sigma$ has dense range, this implies that $\pi\sb
B\circ\rho\sb B=\id\sb B$. For the same reason, $\rho\sb B$ is a
$B$-bimodule morphism. The rest follows from Proposition
\ref{proposition3.7}.
\end{proof}

The following result generalizes \cite[Corollary 1.17]{PiS2007}.

\begin{corollary}
\label{corollary4.15} Let $A$ be a biprojective topological
algebra, and let $I$ be a closed two-sided ideal of $A$ such that
$I$ is complemented as a topological submodule of the
$A$-bi\-mo\-du\-le~$A$. Then $I$ and $A/I$ are biprojective
topological algebras.
\end{corollary}

\begin{proof}
It follows from the hypothesis that there exists another closed
two-sided ideal $J$ of $A$ such that $A=I\oplus J$, and moreover
the operator
\[\lambda\colon I\times J\to A,\quad (a,\,b)\mapsto a+b,\]
is a topological isomorphism.

Since $A$ is biprojective, we have $\overline{A\sp 2}=A$.
This implies that $\overline{A\cdot I}=I$ and
$\overline{A\cdot J}=J$. Hence the topological algebras
\[
I\cong A/J=A/\overline{A\cdot J}\quad\mbox{and}\quad
A/I=A/\overline{A\cdot I}
\]
are biprojective by Proposition \ref{proposition4.14}.
\end{proof}

Now suppose that $A$ and $B$ are topological algebras,
$\varphi\colon A\to B$ is a continuous homomorphism, and
$X$ is a left topological $B$-module. Then, obviously, the
action of $B$ on $X$ induces a continuous action of $A$ on $X$
defined by $a\cdot x=\varphi(a)\cdot x$, where $a\in A$, $x\in X$.
In particular, $B$ is a left topological $A$-module with the
action $a\cdot b=\varphi(a)b$, where $a\in A$, $b\in B$; in this
case, the operator $\varphi\colon A\to B$ is an $A$-module
morphism.

\begin{proposition}\label{proposition4.16}
Let $A$ and $B$ be topological algebras, and let $\varphi\colon
A\to B$ be a continuous homomorphism with dense range. Suppose
further that $A$ is biprojective and has either a left or a right
bounded approximate identity. Then $B$ is biprojective.
\end{proposition}

\begin{proof}
Consider the case where $A$ has a left bounded approximate
identity. Since $A$ is biprojective, there exists an $A$-bimodule
morphism $\rho\colon A\to A\sb +\Ptens A$ such that $\pi\sb
+\circ\rho=\id\sb A$, where ${\pi\sb +\colon A\sb +\Ptens A\to A}$
is the canonical morphism (cf. \cite[Proposition
VII.1.66]{Hel93}). Set $ \rho\sb 1=\id\sb
B\ptens{A}\rho\ptens{A}\id\sb B$. Then
\[\rho\sb 1\colon B\ptens{A} A\ptens{A} B \longrightarrow
B\ptens{A}{}(A\sb + \Ptens A)\ptens{A} B\] is a $B$-bimodule
morphism which is a right inverse for
\[
\pi\sb 1\colon B\Ptens {}(A\ptens{A} B)\to B\ptens{A}
{}(A\ptens{A} B),\quad b\Tens u \mapsto b\Tens\sb Au\qquad (b\in
B,\; u\in A \ptens{A} B).
\]
Since $A$ has a left bounded approximate identity and
$\varphi\colon A\to B$ has dense range, it follows from
Proposition \ref{proposition3.9} that there exists a topological
isomorphism
$$
\kappa\sb B\colon A\ptens{A} B\rightarrow\overline{A\cdot B}=B
$$
uniquely determined by
\[
\kappa\sb B(a\Tens\sb Ab)=a\cdot
b=\varphi(a)b\qquad (a\in A,\; b\in B).
\]
It is obvious that $\kappa\sb B\colon A\ptens{A} B\to B$ is a
right $B$-module morphism. Since $\varphi\colon A\to B$ has dense
range and, clearly (cf. \cite[Chapter XIV, Lemma 1.2]{Mallios}),
$B$ has a left bounded approximate identity, we see that the
$B$-bimodule morphism $\pi\sb 1$ may be identified with the
morphism
\[
\widetilde{\pi}\sb 1\colon B\Ptens B\rightarrow B\ptens{B} B\cong
B,\quad b\Tens c\mapsto bc\qquad (b,c\in B).
\]
But $\pi\sb 1$ (and hence also $\widetilde{\pi}\sb 1$) has a right
inverse $B$-bimodule morphism. So it follows from Proposition
\ref{proposition3.7} that $B$ is biprojective. A similar argument
applies in the case where $A$ has a right bounded approximate
identity.
\end{proof}

The proof of the following proposition repeats the proof of
\cite[Proposition 1.15]{PiS2007}.

\begin{proposition}\label{proposition4.17}
The Cartesian product of an arbitrary family of biprojective
topological algebras is a biprojective topological algebra.
\end{proposition}

The following result generalizes Proposition
\ref{proposition4.17}.

\begin{proposition}
\label{proposition4.18} Let $A=\varprojlim (A\sb \lambda,\,\tau\sp
\mu\sb \lambda\colon A\sb \mu\to A\sb \lambda)$ be a reduced
inverse limit of topological algebras. Suppose that for each
$\lambda$ there exists an $A\sb \lambda$-bimodule morphism
\[\rho\sb \lambda\colon A\sb \lambda\to A\sb \lambda\Ptens A\sb \lambda\] such
that $\pi\sb {A\sb \lambda}\circ\rho\sb \lambda=\id\sb {A\sb
\lambda}$ and that
\[(\tau\sp \mu\sb \lambda\Ptens\tau\sp \mu\sb \lambda)\circ\rho\sb
\mu=\rho\sb \lambda\circ\tau\sp \mu\sb \lambda\]
 whenever $\lambda\prec\mu$. Then $A$ is biprojective.
\end{proposition}

\begin{proof}
Since projective tensor product commutes with reduced inverse
limits (see \cite[\S41.6]{Kothe2}), we may identify $A\Ptens A$
with $\varprojlim (A\sb \lambda\Ptens A\sb \lambda,\,\tau\sp
\mu\sb \lambda\Ptens \tau\sp \mu\sb \lambda)$. It follows from the
assumption that there exists a well-defined continuous linear
operator
\[
\rho=\varprojlim\rho\sb \lambda\colon A\to A\Ptens A.
\]
It is readily verified that $\rho$ is an $A$-bimodule morphism and
that $\pi\sb A\circ\rho=\id\sb A$. Therefore $A$ is biprojective,
as required.
\end{proof}

\begin{ex}\label{example4.19} For
every function $s\colon\N\to\N$, consider the set
\[
M\sb s=\{(i,j)\in\N\times\N \,:\, 1\le j\le s(i)\},
\]
and define
\[
A=\Bigl\{a=(a\sb {ij})\in \CC\sp {\N\times\N} \,:\, \lim\sb
{\substack{(i,j)\to\infty\\ (i,j)\in M\sb s}} |a\sb {ij}|=0,
\quad\forall s\in\N\sp \N\Bigr\}.
\]
Clearly, $A$ is a ${}\sp{\ast}$-subalgebra of $\CC\sp
{\N\times\N}$. Moreover, $A$ is a locally $C\sp{\ast}$-algebra
with respect to the family of $C\sp{\ast}$-seminorms
\[
p\sb s(a)=\sup\sb {(i,j)\in M\sb s} | a\sb {ij}|\qquad (s\in\N\sp
\N)
\]
\emph{(cf.~\cite[Example 5.10]{Phillips1})}. This algebra is
biprojective.
\end{ex}

\begin{proof}
It is easy to see that, for each $s\in\N\sp \N$, the concomitant
$C\sp{\ast}$-algebra $A\sb s=A/\Ker p\sb s$ is isomorphic to $c\sb
0(M\sb s)$, and that the linking maps \[A\sb t\to A\sb s\quad
(s\le t)\] are just the restriction maps from $c\sb 0(M\sb t)$ to
$c\sb 0(M\sb s)$.

For each $(i,j)\in\N\times\N$, let ${e\sb {ij}\in\CC\sp
{\N\times\N}}$ denote the function which is $1$ at $(i,j)$, and
$0$ elsewhere. Fix $s\in\N\sp \N$. By \cite[Example
VII.1.80]{Hel93}, there exists an $A\sb s$-bimodule morphism
$\rho\sb s\colon A\sb s\to A\sb s\Ptens A\sb s$ such that
\[\rho\sb s(e\sb {ij})=e\sb {ij}\Tens e\sb {ij}\] for each $(i,j)\in M\sb s$.
It is readily verified that the family ${\{\rho\sb s\mid s\in\N\sp
\N\}}$ satisfies the conditions of Proposition
\ref{proposition4.18}. This implies that $A$ is biprojective.
\end{proof}

\subsection{Structure theorems for biprojective locally $C\sp{\ast}$-algebras}
Our next aim is to answer the second and the third questions posed
in the beginning of the paper.

\begin{theorem}\label{theorem4.20}
Let $A$ be a locally $C\sp{\ast}$-algebra. Then the following
conditions are equivalent:
\begin{itemize}
\item[\rm(i)] $A$ is a biprojective algebra;
\item[\rm(ii)] all essential left topological $A$-modules are projective;
\item[\rm(iii)] all Hermitian $A$-modules are projective;
\item[\rm(iv)] $\mathrm{Soc}(A)$ is dense in $A$, and moreover,
for each continuous $C\sp{\ast}$-seminorm $p$ on $A$, the
concomitant $C\sp{\ast}$-algebra $A\sb p$ is biprojective;
\item[\rm(v)] $A$ is the direct topological sum of its minimal closed two-sided
ideals each of which is topologically ${}\sp{\ast}$-isomorphic to
a full matrix $C\sp{\ast}$-algebra;
\item[\rm(vi)] $A$ is an annihilator algebra with finite-dimensional
minimal closed two-sided ideals;
\item[\rm(vii)] there exists a
dense ${}\sp{\ast}$-subalgebra $B$ of $A$ which is a biprojective
$C\sp{\ast}$-algebra under some norm and which is continuously
embedded in $A$.
\end{itemize}
\end{theorem}

For the proof of this theorem the following two results will be
important.

\begin{theorem}[{\cite{Hel97,Hel98}}] \label{theorem4.21}
The following properties of a $C\sp{\ast}$-algebra $A$ are
equivalent:
\begin{itemize}
\item[\rm(i)] all Hermitian $A$-modules are projective;
\item[\rm(ii)] $A$ is isometrically ${}\sp{\ast}$-isomorphic to the $c\sb
0$-sum of a family of full matrix $C\sp{\ast}$-algebras.
\end{itemize}
\end{theorem}

The next lemma follows immediately from \cite[Proposition
IV.1.7]{Hel89}.

\begin{lemma}\label{lemma4.22}
Let $A$ and $B$ be locally $C\sp{\ast}$\mbox{-}algebras, and let
$\varphi\colon A\to B$ be a continuous ${}\sp{\ast}$-homomorphism
with dense range. Suppose further that all Hermitian $A$-modules
are projective. Then all Hermitian $B$-modules are projective.
\end{lemma}

\begin{proof}[Proof of Theorem \ref{theorem4.20}] (i)$~\Rightarrow~$(ii)
This follows from Corollary \ref{corollary3.11}.

(ii)$~\Rightarrow~$(iii) This is trivial.

(iii)$~\Rightarrow~$(iv) Since all irreducible Hermitian
$A$-modules are projective, it follows from Theorem
\ref{theorem3.24} that $\mathrm{Soc}(A)$ is dense in $A$. Further,
let $p$ be a continuous $C\sp{\ast}$-seminorm on $A$, and let
$A\sb p$ be the corresponding concomitant $C\sp{\ast}$-algebra.
Using Lemma \ref{lemma4.22}, we see that all Hermitian $A\sb
p$-modules are projective.  Applying Theorem \ref{theorem4.21} and
taking into account Example \ref{example4.4}, we conclude that
$A\sb p$ is biprojective.

(iv)$~\Rightarrow~$(v) This follows from Theorem \ref{theorem3.3}
(see also Theorem \ref{theorem3.24}) and Theorem \ref{theorem4.5}.

(v)$~\Rightarrow~$(vi) This follows from Theorem
\ref{theorem3.24}.

(vi)$~\Rightarrow~$(vii) This follows from Theorem
\ref{theorem3.24}, Theorem \ref{theorem3.3}(xi), and
Example \ref{example4.4}.

(vii)$~\Rightarrow~$(i) This follows from Proposition
\ref{proposition4.16}.
\end{proof}

\begin{proposition}
\label{proposition4.23} Let $A$ be a biprojective locally
$C\sp{\ast}$-algebra. Then the following statements hold.
\begin{enumerate}
\item[\rm(i)] For each closed left ideal $J$ of
$A$, $J$ and $A/J$ are projective left topological $A$-modules.
Moreover, $J$ is complemented as a topological submodule of the left $A$-module $A$.
\item[\rm(ii)] For each closed right ideal $K$ of $A$, $K$ and $A/K$ are
projective right topological $A$-modules. Moreover, $K$ is complemented as a
topological submodule of the right $A$-module $A$.
\item[\rm(iii)]
For each closed two-sided ideal $I$ of $A$, $I$ and $A/I$ are
projective topological $A$\mbox{-}bi\-mo\-du\-les. Moreover, $I$
is complemented as a
topological submodule of the $A$-bimodule~$A$.
\end{enumerate}
\end{proposition}

\begin{proof} (i) It follows from Theorem \ref{theorem4.20} that $A$ is an annihilator
algebra, and that $\mathrm{Soc}(A)$ is dense in $A$.  Now Theorem
\ref{theorem3.24} implies that $J$ is complemented as a topological
submodule of the left $A$-module $A$. Hence $J$ and
$A/J$ are retracts (i.e., $A$-module direct summands) of the module $A$.
Since the biprojective algebra $A$ is left projective, we see that
$J$ and $A/J$ are projective left topological $A$-modules (see
\cite[Proposition III.1.16]{Hel89}).

(ii) and (iii) These are similar.
\end{proof}

The next result follows easily from Theorems \ref{theorem3.3} and
\ref{theorem4.20} (see also Proposition \ref{proposition4.23} and
Corollary \ref{corollary4.15}).

\begin{proposition}\label{proposition4.24}
Let $A$ be a biprojective locally $C\sp{\ast}$-algebra, and let
$I$ be a closed two-sided ideal of $A$. Then $I$ and $A/I$ are
biprojective locally $C\sp{\ast}$-algebras, and moreover $A$ is
topologically ${}\sp{\ast}$-isomorphic to the Cartesian product of
the algebras $I$ and $A/I$.
\end{proposition}

Propositions \ref{proposition3.1} and \ref{proposition4.24} (see
also Theorem \ref{theorem4.20}) imply the following.

\begin{corollary}
Let $A$ be a biprojective locally $C\sp{\ast}$-algebra, and let
$P$ be the family of all continuous $C\sp{\ast}$-seminorms on $A$.
Then, for each $p\in P$, the concomitant $C\sp{\ast}$-algebra
$A\sb p$ is biprojective and, as a consequence, $A$ can be
represented as an inverse limit of biprojective
$C\sp{\ast}$-algebras.
\end{corollary}

Applying Theorems \ref{theorem4.20} and \ref{theorem3.29} and
taking into account Example \ref{example4.1}, we get the following
theorem.

\begin{theorem}\label{theorem4.26}
Let $A$ be a unital locally $C\sp{\ast}$-algebra. Then the
following conditions are equivalent:
\begin{itemize}
\item[\rm(i)] $A$ is a biprojective (or, equivalently, contractible) algebra;
\item[\rm(ii)] all essential (or, equivalently%
\footnote[16]{\thinspace Clearly, every essential left topological
module over a
unital topological algebra is unital.}%
, unital) left topological $A$-modules are projective;
\item[\rm(iii)] all Hermitian $A$-modules are projective;
\item[\rm(iv)] all irreducible Hermitian $A$-modules are projective;
\item[\rm(v)] $\mathrm{Soc}(A)$ is dense in $A$;
\item[\rm(vi)] $A$ is an annihilator algebra;
\item[\rm(vii)] $A$ is topologically ${}\sp{\ast}$-isomorphic to the
Cartesian product of a family of full matrix algebras.
\end{itemize}
\end{theorem}

\begin{remark}\label{remark4.27}
The equivalence of conditions (i), (ii) and (vii) of Theorem
\ref{theorem4.26} was proved earlier in \cite[Theorem 3.3]{Frag2}.
\end{remark}

In \cite[Corollary 5.4]{PiS2007}, we proved that every
biprojective $\sigma$-$C\sp{\ast}$-algebra is topologically
${}\sp{\ast}$\mbox{-}iso\-morphic to the Cartesian product of a
family of biprojective $C\sp{\ast}$-algebras. In the same paper,
we have raised the following question (see
\cite[Question~5.5]{PiS2007}):

\medskip

\emph{Can this result be extended to arbitrary (i.e.,
non-metrizable) locally $C\sp{\ast}$-algebras?}

\medskip

Now we show that the answer to this question is negative. Namely,
we present an example of a non-unital biprojective locally
$C\sp{\ast}$-algebra that is not topologically isomorphic to a
Cartesian product of (biprojective) $C\sp{\ast}$-algebras.

\begin{theorem}\label{theorem4.28}
There exists a biprojective locally $C\sp{\ast}$-algebra that is
not topologically isomorphic to a Cartesian product of
$C\sp{\ast}$-algebras.
\end{theorem}

\begin{proof}
Consider the biprojective locally $C\sp{\ast}$-algebra $A$ defined
in Example \ref{example4.19}. Assume towards a contradiction that
this algebra is topologically isomorphic to the Cartesian product
$\prod\limits\sb {i\in I} A\sb i$ of a family of
$C\sp{\ast}$-algebras $A\sb i$, $i\in I$. Since $A$ is
biprojective, it follows from Proposition \ref{proposition4.24}
that $A\sb i$ is biprojective for each $i\in I$. Since $A$ is
commutative, it follows from Theorem \ref{theorem4.5} that each
algebra $A\sb i$, $i\in I$, is isometrically
${}\sp{\ast}$-isomorphic to a $C\sp{\ast}$-algebra of the form
$c\sb 0(\Lambda\sb i)$ for some set $\Lambda\sb i$. So $A$ is
topologically isomorphic to the algebra $\prod\limits\sb {i\in I}
c\sb 0(\Lambda\sb i)$.

Denote the latter algebra by $B$. Note that $A$ is not metrizable,
because the set $\N\sp \N$ of all functions $s\colon\N\to\N$ does
not have a countable cofinal subset (cf.~\cite[Example
5.10]{Phillips1}). Therefore $B$ is not metrizable either, and so
$I$ is uncountable. For each $i\in I$, let \[p\sb i\colon B\to
c\sb 0(\Lambda\sb i)\] denote the natural projection onto the
$i$th factor. Then it is easy to check (cf.~\cite[Proposition
V.1.8]{Hel93}) that every non-zero continuous character on $B$ has
the form \[b\mapsto p\sb i(b)(\nu)\] for some $i\in I$ and for
some $\nu\in \Lambda\sb i$. Therefore there is no countable
separating set of continuous characters on $B$. On the other hand,
the characters
\[
\{a\mapsto a\sb {ij} \mid (i,j)\in\N\times\N\}
\]
obviously separate the points of $A$. The resulting contradiction
shows that $A$ is not topologically isomorphic to $B$. This
completes the proof.
\end{proof}

Thus Theorem \ref{theorem4.7} and \cite[Corollary 5.4]{PiS2007}
can not be generalized to arbitrary locally $C\sp{\ast}$-algebras.
Nevertheless we have the following result.

\begin{proposition}
Let $A$ be a non-unital biprojective locally $C\sp{\ast}$-algebra.
Then $A$ is topologically ${}\sp{\ast}$-isomorphic to the
Cartesian product of two biprojective locally
$C\sp{\ast}$-algebras one~of which is an infinite-dimensional
$C\sp{\ast}$-algebra.
\end{proposition}

\begin{proof}
By Theorem \ref{theorem4.20}, $A$ is an annihilator algebra. Now
Proposition \ref{proposition3.33} implies that $A$ is
topologically ${}\sp{\ast}$-isomorphic to the Cartesian product of
two annihilator locally $C\sp{\ast}$\mbox{-}algebras one of which
is an infinite-dimensional $C\sp{\ast}$-algebra. The rest follows
from Proposition \ref{proposition4.24}.
\end{proof}

\section{Superbiprojective locally $C\sp{\ast}$-algebras}
Here our aim is to answer the fourth question posed in the
beginning of the paper.

The following result is well known (cf. \cite[Corollary
8]{Sel96C}, see also \cite[Corollary 2.4.3]{Sel2002}).

\begin{theorem}\label{theorem5.1}
Let $A$ be a biprojective topological algebra. Then the following
conditions are equivalent:
\begin{enumerate}

\item[\rm(i)]
$\mathcal{H}\sp 2(A,X)=0$ for all topological $A$-bimodules $X$;
\item[\rm(ii)]
$\mathcal{H}\sp 2(A,A\Ptens A)=0$;
\item[\rm(ii)]
the operator
\[
\Delta\colon A\Ptens A\longrightarrow (A\sb +\Ptens A)\oplus {}
(A\Ptens A\sb +),\quad a\Tens b\mapsto (a\Tens b,\,a\Tens b),
\]
has a left inverse $A$-bimodule morphism.

\end{enumerate}
\end{theorem}

\begin{definition}[{cf.\ \cite{Sel2001} and \cite{PiS2007}}]
A topological algebra $A$ is said to be \emph{superbiprojective}
if $A$ is biprojective and satisfies the equivalent conditions of
Theorem~\ref{theorem5.1}.
\end{definition}

It is clear that every contractible topological algebra is
superbiprojective.

\medskip

The proof of the next proposition repeats, in essence, the proof
of \cite[Theorem 1]{Sel75}.

\begin{proposition}
Let $A$ be a biprojective topological algebra. If $A$ has a left
or a right identity, then $A$ is superbiprojective.
\end{proposition}

\begin{example}
Let $F$ be a complete Hausdorff locally convex space with $\dim
F>1$. Then every continuous linear functional $f\in F\sp{\ast}$,
$f\neq 0$, defines on $F$ the structure of a topological algebra,
denoted by $S\sb f(F)$, with multiplication given by $ab=f(a)b$.
It is easily seen that $S\sb f(F)$ is a superbiprojective algebra
with a left identity (\cite[Example~3]{Sel96B}, see also
\cite[Theorem 1]{Sel75} and \cite[Example 2.6]{Pir2004}). Clearly,
$S\sb f(F)$ is not contractible. Note that $S\sb f(F)$ is a
particular case of the algebras $E\Ptens F$ considered in Example
\ref{example4.11}. Indeed, it suffices to set $E=\CC$ and to
define $\langle \lambda ,x\rangle = \lambda f(x)$ for
$\lambda\in\CC$ and $x\in F$.
\end{example}

The following result generalizes \cite[Proposition 2.8]{PiS2007}.

\begin{proposition}
Let $A$ be a superbiprojective topological algebra, and let $I$ be
a closed two-sided ideal of $A$ such that $I$ is complemented as a
topological submodule of the $A$\mbox{-}bimodule $A$. Then $I$ and
$A/I$ are superbiprojective topological algebras.
\end{proposition}

\begin{proof}
Using Corollary \ref{corollary4.15}, we see that $I$ and $A/I$ are
biprojective topological algebras. Also, it follows from the
hypothesis that there is a closed two-sided ideal $J$ of $A$ such
that $A=I\oplus J$, and the operator $\lambda\colon I\times J\to
A$, $(a,\,b)\mapsto a+b$, is a topological isomorphism. Moreover,
as we noted in the proof of Corollary \ref{corollary4.15},
$\overline{A\cdot I}=I$ and $\overline{A\cdot J}=J$. Now
Proposition \ref{proposition3.13} implies that the left
topological $A$-modules
\[ I\cong A/J=A/\overline{A\cdot
J}\quad\mbox{and}\quad A/I=A/\overline{A\cdot I}
\]
are projective. Similarly, $I\cong A/J$ and $A/I$ are projective
in the category of right topological $A$-modules. Therefore the
natural projections of topological algebras $A\to A/J$ and $A\to
A/I$ satisfy the conditions of \cite[Theorem 5.2]{Pir2002}, and so
the algebras $I\cong A/J$ and $A/I$ are superbiprojective.
\end{proof}

By combining this result with Proposition
\ref{proposition4.23}(iii) we get the following.

\begin{corollary}
Let $A$ be a superbiprojective locally $C\sp{\ast}$-algebra, and
let $I$ be a closed two-sided ideal of $A$. Then $I$ and $A/I$ are
superbiprojective locally $C\sp{\ast}$-algebras.
\end{corollary}

This result, together with Proposition
\ref{proposition3.1}, implies the following.

\begin{corollary}\label{corollary5.7} Let $A$ be a superbiprojective
locally $C\sp{\ast}$-algebra, and let $P$ be the family of all
continuous $C\sp{\ast}$-seminorms on $A$. Then, for each $p\in P$,
the concomitant $C\sp{\ast}$-algebra $A\sb p$ is
superbiprojective.
\end{corollary}

Now we recall a result concerning $C\sp{\ast}$-algebras.

\begin{theorem}[{see \cite[Corollary 4.64]{Sel2000} and
\cite{Sel2001}}]\label{theorem5.9} Let $A$ be a superbiprojective
$C\sp{\ast}$-algebra. Then $A$ is topologically
${}\sp{\ast}$-isomorphic to the Cartesian product of finitely many
full matrix algebras. In particular, $A$ is finite-dimensional and
contractible.
\end{theorem}

Everything is now ready for the proof of the main result of this
section. The following theorem generalizes \cite[Theorem 5.6 and
Corollary 5.7]{PiS2007}.

\begin{theorem}\label{theorem5.10}
Let $A$ be a locally $C\sp{\ast}$-algebra. Then the following
conditions are equivalent:
\begin{itemize}
\item[\rm(i)] $A$ is a superbiprojective algebra;
\item[\rm(ii)] $A$ is topologically ${}\sp{\ast}$-isomorphic to
the Cartesian product of a family of full matrix algebras;
\item[\rm(iii)] $A$ is a contractible algebra.
\end{itemize}
\end{theorem}

\begin{proof} (i)$~\Rightarrow~$(ii) It follows from the hypothesis
and from Theorem \ref{theorem4.20} that $A$ is biprojective, and
that $\mathrm{Soc}(A)$ is dense in $A$. On the other hand,
Corollary \ref{corollary5.7} and Theorem~\ref{theorem5.9} imply
that, for each continuous $C\sp{\ast}$-seminorm $p$ on $A$, the
concomitant $C\sp{\ast}$-algebra $A\sb p$ is finite-dimensional.
Now the result follows from Corollary \ref{corollary3.4}(ii).

(ii)$~\Rightarrow~$(iii) This follows from \cite[Lemma
11]{Sel96A} and
\cite[Proposition~VII.1.73]{Hel93}. (See also Example
\ref{example4.1}.)

(iii)$~\Rightarrow~$(i) This is trivial.
\end{proof}

Recall that the equivalence of conditions (ii) and (iii) in
Theorem \ref{theorem5.10} was proved earlier in \cite[Theorem
3.3]{Frag2} (see Remark \ref{remark4.27}).

\medskip

We now recall an important definition (see \cite{Hel89,Hel2000}).

\begin{definition}
Let $A$ be a topological algebra. The least $n$ for which
$\mathcal{H}\sp m(A,X)=0$ for all topological $A$-bimodules $X$
and all $m>n$, or $\infty$ if there is no such $n$, is called the
\emph{homological bidimension} of $A$, and is denoted by
$\mathop{\mathrm{db}}A$.
\end{definition}

Clearly, $\mathop{\mathrm{db}}A=0$ means exactly that $A$ is
contractible. Also, a biprojective topological algebra $A$ always
has $\mathop{\mathrm{db}}A\le 2$, and it has
$\mathop{\mathrm{db}}A\le 1$ if and only if it is
superbiprojective.

\begin{example}
If $G$ is an infinite compact Lie group, then the Fr\'echet
algebra $\mathcal{E}(G)$ (see Example~\ref{example4.9}) is
superbiprojective and $\mathop{\mathrm{db}}\mathcal E(G)=1$
\cite[Corollary 5.4]{Pir2002}.
\end{example}

From Theorem \ref{theorem5.10} we immediately get the following
corollaries.

\begin{corollary}
Let $A$ be a non-unital biprojective locally $C\sp{\ast}$-algebra.
Then $\mathop{\rm db}A=2$ and $\mathcal{H}\sp 2(A,A\Ptens A)\neq
0$.
\end{corollary}

\begin{corollary}
In the class of biprojective locally $C\sp{\ast}$-algebras the
homological bidimension $\mathop{\mathrm{db}}A$ can take only the
values\/ $0$ and\/ $2$.
\end{corollary}

\section*{Acknowledgements}

This paper is based on a lecture of the second author delivered at
the 19$\sp \mathrm{th}$ International Conference on Banach
Algebras held at
{\fontencoding{T1}\selectfont B\k{e}dlewo}, July 14--24, 2009. The support for
the meeting by the Polish Academy of Sciences, the European
Science Foundation under the ESF-EMS-ERCOM partnership, and the
Faculty of Mathematics and Computer Science of the Adam Mickiewicz
University at Pozna\'n is gratefully acknowledged.

The research was partially supported by the RFBR grant
08-01-00867. The first author was also partially supported by the
Ministry of Education and Science of Russia (programme
``Development of the scientific potential of the Higher School'',
grant no. 2.1.1/2775), and by the President of Russia grant
MK-1173.2009.1. We are grateful to D.~P.~Blecher,
M.~Fragoulopoulou, M.~Haralampidou, A.~Ya.~Helemskii, and
G.~Racher for useful discussions. Also, we are grateful to the
referee for his/her valuable comments.

\end{document}